\documentclass[11pt]{amsart}
\usepackage{latexsym}
\usepackage{amssymb,amsmath,amsthm,amscd, amssymb,enumerate, verbatim}
\def\co{\colon\thinspace}

\newcommand{\vol}{\mbox{vol}}

\newcommand{\Int}{\mbox{Int}}

\newcommand{\e}{\epsilon}
\def\d{\partial}
\def\a{\alpha}
\def\G{\Gamma}

\def\H{\mathcal H}
\def\G{\mathcal G}
\def\S{\mathcal S}

\def\i{\iota}
\def\e{\epsilon}

\def\s{\sigma}

\def\R{\mathbb{R}}
\def\Z{\mathbb{Z}}
\newcommand{\g}{\gamma}

\newtheorem{thm}{Theorem}[section]
\newtheorem{cor}[thm]{Corollary}
\newtheorem{quest}[thm]{Question}
\newtheorem{lem}[thm]{Lemma}
\newtheorem{prop}[thm]{Proposition}

\newtheorem{assum}[thm]{Assumption}
\newtheorem{Example}[thm]{Example}
\newenvironment{ex}{\begin{Example}\rm}{\end{Example}}
\newtheorem{remark}[thm]{Remark}
\newenvironment{rmk}{\begin{remark}\rm}{\end{remark}}
\newtheorem{Fact}[thm]{Fact}

\newtheorem{Main Lemma}[thm]{Main Lemma}

\newtheorem{Convention}[thm]{Convention}

\begin{document}\abovedisplayskip=6pt plus3pt minus3pt 
\belowdisplayskip=6pt plus3pt minus3pt
\title[Aspherical manifolds and
relatively hyperbolic groups] 
{\bf Aspherical manifolds
with relatively hyperbolic fundamental groups\rm}
\date{}
\rm

\author{Igor Belegradek}
\address{Igor Belegradek\\School of Mathematics\\ Georgia Institute of
Technology\\ Atlanta, GA 30332-0160}\email{ib@math.gatech.edu}

\begin{abstract} 
We show that the aspherical manifolds produced via the relative strict 
hyperbolization of polyhedra enjoy many group-theoretic and topological   
properties of open finite volume negatively pinched manifolds, including 
relative hyperbolicity, nonvanishing of simplicial volume, co-Hopf property, 
finiteness of outer automorphism group, absence of splitting over elementary 
subgroups, and acylindricity. 
In fact, some of these properties hold for any compact aspherical 
manifold with incompressible aspherical boundary
components, provided the fundamental group is hyperbolic relative to
fundamental groups of boundary components. 
We also show that no manifold obtained via the relative strict 
hyperbolization can be embedded into a compact
K\"ahler manifold of the same dimension, except when the
dimension is two.
\end{abstract} 
\maketitle

\section{Introduction}
\label{sec: intro}

The subject of negatively curved manifolds sits at the crossroads
of Riemannian geometry, geometric group theory, coarse topology,
and dynamics, and there is an extensive literature 
on (complete Riemannian) negatively curved 
manifolds of finite volume. 
In dimensions $>3$ the main source of examples is
arithmetic lattices; 
other known examples are typically obtained from arithmetic ones
by various ``surgery'' constructions along totally geodesic
submanifolds~\cite{P-SGro, MosSiu, GroThu, Deraux, FarJon, FJO}. 
These constructions are very special,
and one could instead try to work in a broader context of 
aspherical manifolds whose fundamental groups are hyperbolic or 
relatively hyperbolic hoping, in particular, that this could 
eventually lead to other sources of examples.

The strict hyperbolization of Charney-Davis~\cite{CD} gives a rich
supply of closed aspherical manifolds with hyperbolic fundamental groups.
This is a 
functorial procedure that turns any finite simplicial complex $K$
into a locally-$CAT(-1)$ piecewise  hyperbolic finite cell 
complex $\H(K)$. 
A key property of the procedure is that it preserves links of all
the faces, and hence if $K$ is a triangulated
manifold, then so is $\H(K)$ (and this also works in smooth 
category). 
The closed locally-$CAT(-1)$ piecewise hyperbolic manifolds 
obtained via the strict hyperbolization behave in many 
ways like closed Riemannian manifolds of negative sectional 
curvature; in fact, if $K$ is a smoothable
triangulation of a closed manifold, then for all we know,
$\H(K)$ could admit a negatively curved Riemannian metric,
and this is what happens e.g. in dimensions $\le 3$ 
(by Thurston's hyperbolization theorem). 

The strict hyperbolization has a relative 
version~\cite{DJW, CD} (cf.~\cite{DJ, Gro-hgr}), 
and the main theme of this paper is establishing 
an analogy between the aspherical manifolds produced via
relative strict hyperbolization, and the open complete finite
volume negatively pinched manifolds. 
Roughly, the relative strict hyperbolization,
which is explained in detail in 
Section~\ref{sec: strict rel hyper} 
is a procedure that takes as the input a pair of finite
simplicial complexes $(K,L)$, where $L$ is a nonempty
subcomplex of $K$,
and produces a pair of simplicial complexes $(R_K, R_L)$
with the following properties: 
\begin{itemize}
\item $R_L$ is isomorphic to a subdivision of $L$,
\item
$R_K$ is aspherical if and only if 
each path-component of $L$ is aspherical;
\item
$R_L$ is incompressible
in $R_K$ (i.e. no homotopically nontrivial
loop in $R_L$ is null-homotopic in $R_K$);
\item
if $K$ is a (triangulated, smooth or PL)
manifold with boundary $L$, then $R_K$ is a 
(respectively, triangulated, smooth or PL) 
manifold with boundary $R_L$. 
\item
the space $Z$ obtained from $R_K$ by attaching a cone 
to each path-component of $R_L$ admits a piecewise-hyperbolic
locally-$CAT(-1)$ metric.
\end{itemize}

An appealing feature of relative strict hyperbolization is that
it allows to construct compact aspherical manifolds with prescribed
boundary: namely, if $L$ bounds and each path-component of $L$ is 
aspherical, then $L$ bounds an aspherical manifold. By contrast,
the problem of prescribing topology of cusps of negatively pinched
finite volume complete Riemannian manifolds is still wide open; 
for recent progress in locally-symmetric case 
see~\cite{LonRei-eta, LonRei-orbi, McR}.

We emphasize that $R_K$ itself carries no obvious 
locally-$CAT(-1)$ metric, in fact, geometrically,
$R_K$ looks more like a finite volume negatively pinched manifold
with cusps chopped off, with $R_L$ corresponding to
the union of cusp cross-sections on the boundary. 
However, negative curvature
persists on the group-theoretic level as follows.

\begin{thm}\label{intro-thm: rel hyp}
If $(R_K, R_L)$ is a pair of finite (nonempty) simplicial complexes
obtained via relative strict hyperbolization, and $L_1,\dots, L_m$
are all the path-components of $R_L$ that have infinite fundamental
groups, then $\pi_1(R_K)$ is a non-elementary relatively hyperbolic 
group in the sense of Bowditch,
specifically, $\pi_1(R_K)$ is hyperbolic relative to 
the subgroups $\pi_1(L_1),\dots, \pi_1(L_m)$. In particular, 
if each component of $R_L$ has finite fundamental group, 
then $\pi_1(R_K)$ is a non-elementary hyperbolic group. 
\end{thm}

Relatively hyperbolic groups 
(see Section~\ref{sec: proving rel hyp} for a definition)
were introduced by Gromov~\cite{Gro-hgr}
and since then various characterizations of relatively
hyperbolic groups have been obtained by Farb~\cite{Far-rel}, 
Bowditch~\cite{Bow-rel}, Yaman~\cite{Yaman}, 
Osin~\cite{Osi-rel}, Dru{\c{t}}u-Osin-Sapir~\cite{DOS}. 
Recently there has been a surge of activity in this subject,
and many results about hyperbolic groups have been adapted to
the relatively hyperbolic setting.  
Three basic classes of examples of relatively hyperbolic groups 
are free products with finitely many factors 
(which are hyperbolic relative to the factors), hyperbolic groups
(which are hyperbolic relative to the empty set of subgroups),
and geometrically finite isometry groups of negatively
pinched Hadamard manifolds (which are hyperbolic relative
to the maximal parabolic subgroups). Finally,
any group is hyperbolic relative to itself. A subgroup of
a relatively hyperbolic group is called {\it elementary} if it is 
finite, or virtually-$\Z$, or hyperbolic relative 
to itself; otheriwise, the subgroup is called {\it non-elementary}. 
Another major source of examples of relatively 
hyperbolic groups is the small cancellation 
theory~\cite{Osi-sc}.
There are however very few (primary) constructions of higher-dimensional 
relatively hyperbolic groups, e.g. the fundamental groups of 
compact aspherical manifolds that are hyperbolic
relative to the fundamental groups of the boundary components,
and this is exactly what Theorem~\ref{intro-thm: rel hyp} does.

Theorem~\ref{intro-thm: rel hyp} has a complicated history.
It was predicted by Gromov, and first stated (without proof) 
in~\cite[p. 257]{Gro-hgr} for a different hyperbolization 
procedure. However, results of Charney-Davis~\cite{CD} imply
that the claim in~\cite[p. 257]{Gro-hgr} was 
incorrect because the Gromov's hyperbolization procedure 
does not yield negative curvature in dimensions $\ge 4$.
This motivated Charney-Davis to introduce the strict
hyperbolization whose relative version we study in this paper.
Theorem~\ref{intro-thm: rel hyp} answers a question of
Szczepa{\'n}ski~\cite{Szc-rel} who proved a weaker 
result that $\pi_1(R_K)$ of Theorem~\ref{intro-thm: rel hyp} is 
relatively hyperbolic in the sense of Farb, 
relative to the same collection of subgroups.
Goldfarb in~\cite[Example 5.5]{Gol-nov} claimed a result 
that is stronger than Theorem~\ref{intro-thm: rel hyp}
(in the special case when $R_K$ is an aspherical manifold
with boundary $R_L$), yet he recently acknowledged that
his proof is incorrect; I do not know how to 
prove what Goldfarb claimed. 

We next turn to a general study of aspherical manifolds
with relatively hyperbolic fundamental groups; more precisely
we make the following assumption.

\begin{assum}\label{intro: assum} 
For $n>2$, $M$ denotes a compact aspherical $n$-manifold 
with non-empty boundary $\d M$
such that each component of $\d M$ is aspherical and incompressible
in $M$, and $\pi_1(M)$ is non-elementary relatively hyperbolic, 
relative to the fundamental groups of the components of $\d M$. 
\end{assum}

Here are examples of manifolds that {\it do} or {\it do not} satisfy
Assumption~\ref{intro: assum}.
\begin{itemize}
\item
If $R_K$ is a compact $n$-manifold with boundary $R_L$, and
if each component of $R_L$ is aspherical, 
then by Theorem~\ref{intro-thm: rel hyp}, 
$R_K$ satisfies Assumption~\ref{intro: assum}.
\item
Any finite volume complete negatively pinched (e.g. of constant negative
curvature) Riemannian manifold with totally geodesic boundary
satisfies Assumption~\ref{intro: assum}; 
as we explain in Section~\ref{sec: geod bound}
this essentially follows from~\cite[Section 7]{Bow-rel}.
\item 
If $M$ is the quotient of a nonpositively curved symmetric space
of rank $\ge 2$ by a non-cocompact torsion-free lattice, then 
$M$ does not satisfy Assumption~\ref{intro: assum}~\cite{BerDruMos}.
\item
If $M$ is a $3$-manifold that satisfies Assumption~\ref{intro: assum},
then we show in Section~\ref{sec: acyl for maps} that $M$ is acylindrical. 
Note that there are non-acylindrical $3$-manifold with 
incompressible boundary whose interior admits a complete
hyperbolic metric (see e.g.~\cite[Example 1.4.5]{CanMcC}).
\end{itemize}

If $G$ is the fundamental group of a 
closed aspherical manifold, and if $G$ is hyperbolic,
then $G$ admits no nontrivial splitting
over finite, or virtually-$\Z$ subgroups (as immediately follows 
from the Mayer-Vietoris sequence in group cohomology), so by
Rips theory~\cite{BesFei} 
$G$ admit no nontrivial isometric actions on 
$\R$-trees with finite, or virtually-$\Z$ arc stabilizers, and
in particular, $G$ is co-Hopfian~\cite{RipSel}
and $\mathrm{Out}(G)$ is finite~\cite{Pau}. 

To state a
relatively hyperbolic version of these results we
need the following definition.
We say that a group $G$ {\it has property} 
(m) if $G$ is hyperbolic relative to a family of 
subgroups none of which is isomorphic
to a non-elementary relatively hyperbolic group.
For example, if 
$G$ is not isomorphic to a non-elementary 
relatively hyperbolic group, then $G$ has property (m),
where we consider $G$ hyperbolic relative to itself.
It was proved in~\cite{BerDruMos} that Dunwoody's inaccessible 
group does not have property (m).
It is unknown whether there exists a torsion-free
or finitely presented group 
that does not have property (m), and it is conceivable that 
the fundamental group of a closed aspherical manifold
always has property (m).

\begin{thm} \label{intro-thm: no split}
If $M$ satisfies Assumption~\ref{intro: assum}, then \newline
\rm (i)\it\ 
$\pi_1(M)$ admits no nontrivial splitting
over an elementary 
subgroup.\newline
\rm (ii)\it\  
$\pi_1(M)$ is co-Hopfian.\newline
\rm (iii)\it\ 
If the fundamental group of every component of $\d M$
has property \textup{(m)}, then
$\mathrm{Out}(\pi_1(M))$ is finite. 
\end{thm}

Part (i) of Theorem~\ref{intro-thm: no split}
is proved as in~\cite{Bel-Top} where we obtained
a similar result for finite volume negatively pinched manifolds,
while parts (ii)-(iii)
easily follow from (i) and the recent remarkable work of 
Dru{\c{t}}u-Sapir~\cite{DruSap-rips} that extends Rips theory
to the relatively hyperbolic setting. 

Here is an open question on automorphisms of $\pi_1(M)$.

\begin{quest}\label{quest: rh aut pres parab}
Let $M$ satisfy Assumption~\ref{intro: assum}. 
Does every automorphism of $\pi_1(M)$ permutes the maximal parabolic
subgroups?
\end{quest}

Another approach to proving part (ii) of 
Theorem~\ref{intro-thm: no split} is based on recent work
of Mineyev-Yaman~\cite{MinYam}, who
showed among other things that if $M$ satisfies 
Assumption~\ref{intro: assum}, then $||M,\d M||>0$ where
$||M,\d M||$ is the relative simplicial volume of $(M,\d M)$. 
We use this result of Mineyev-Yaman to prove the following. 

\begin{thm}\label{thm: finite coker}
If $M_1$, $M_2$ are $n$-manifolds that
satisfy Assumption~\ref{intro: assum}, 
and $\phi\co\pi_1(M_1)\to\pi_1(M_2)$ is an injective
homomorphism that maps each maximal parabolic subgroup 
to a parabolic subgroup, then $\phi(\pi_1(R_1))$ has 
finite index in $\pi_1(R_2)$. 
If in addition $M_1$ $M_2$ are homeomorphic, 
then $\phi$ is an isomorphism.
\end{thm}

We call a group $\G$ {\it intrinsically elementary} if the image of
any injective homomorphism of $\G$ into a relatively hyperbolic group
is elementary (i.e. finite, virtually-$\Z$, or parabolic).
Thus if in Theorem~\ref{thm: finite coker} all maximal
parabolic subgroups of $\pi_1(M_1)$ are intrinsically elementary,
and $n>2$, then the assumption ``$\phi$ maps each maximal parabolic subgroup 
to a parabolic subgroup'' holds true. 
(To see that in this case ``elementary'' implies ``parabolic'', note 
that each boundary component of $M_1$ is a closed aspherical manifold
of dimension $\ge 2$ so its fundamental group 
cannot be finite or virtually-$\Z$). 

Below we list several classes of closed aspherical manifolds
with intrinsically elementary fundamental groups, and these 
manifolds can be realized (via the relative strict hyperbolization)
as boundary components of manifolds satisfying 
Assumption~\ref{intro: assum}.
Examples of intrinsically elementary groups include
\begin{itemize}
\item
groups with no nonabelian free subgroups~\cite{Tuk}
(e.g. fundamental groups of infrasolvemanifolds), 
\item
fundamental groups of closed $3$-dimensional graph 
manifolds~\cite{BerDruMos},
\item
groups with finite
dimensional second bounded cohomology~\cite{Fuj} 
(e.g. cocompact 
irreducible lattices in higher
rank Lie groups~\cite{Burger-Monod, Burger-Monod-err}).
\item
Since any infinite elementary group has elementary normalizer,
any group that
contains an intrinsically elementary infinite normal subgroup
is intrinsically elementary. For example,
the fundamental group of any fiber bundle with aspherical base 
and fiber is intrinsically elementary, provided the fiber has
intrinsically elementary fundamental groups, and this construction
can be iterated since the total space of any such bundle is aspherical. 
\end{itemize}

We next focus on (not necessarily aspherical) 
compact $n$-manifolds obtained
via relative strict hyperbolization. Thus $K$ is an $n$-manifold
with boundary $L$, so that $R_K$ is 
an $n$-manifold with boundary $R_L$. 
To emphasize that we are dealing with manifolds,
we denote $(R_K,R_L)$ by $(R,\d R)$ and say that {\it $R$ is 
an $n$-manifold obtained via relative strict hyperbolization}.
We give $\pi_1(R)$ the structure of a relatively hyperbolic
group given by Theorem~\ref{intro-thm: rel hyp}, so that
the conjugacy classes of maximal parabolic subgroups bijectively
correspond to the components of $\d R$ that have infinite 
fundamental group. 

According to~\cite{Grom-vol-bounded-coh},
there exist positive constants $C_1=C_1(n)$, 
$C_2=C_2(n,a)$ such that if $N$ is a compact $n$-manifold and $V=\Int(N)$
admits a complete finite volume Riemannian metric of sectional
curvature within $[-a^2, -1]$, then
$C_1\vol(V)\ge ||N,\d N||\ge C_2\vol(V)$. We prove a similar
statement for manifolds obtained via the
relative strict hyperbolization.

\begin{thm} \label{intro: thm-hausd}
There are constants $C_1\ge C_2>0$ depending
only on $n$ such that
if $R$ is a compact $n$-manifold obtained by the
relative strict hyperbolization, then
$C_1\mathcal H^n(R)\ge ||R,\d R||\ge C_2\mathcal H^n(R)$. 
\end{thm}

Here $\mathcal H^n(R)$ denotes 
the $n$-dimensional Hausdorff measure of $R$, and the metric
on $R$ is induced by the inclusion $R\to Z$.
As we explain in Section~\ref{sec: hyp and simp vol}, 
$\mathcal H^n(R)$ is roughly the same as the number of 
$n$-simplices in some natural triangulation of $R$.
Since there are only finitely many ways to glue
finitely many simplices, we get the following finiteness
theorem (which also holds in PL and smooth categories).

\begin{thm}\label{intro-thm: finiteness}
For any $C>0$, 
there are only finitely many homeomorphism
types of triangulated $n$-manifolds $R$ obtained by the relative strict
hyperbolization and satisfying satisfying $||R,\d R||<C$.
\end{thm}

One immediate implication of Theorem~\ref{intro: thm-hausd} 
is that $||R,\d R||>0$. (If $R$ is aspherical, this also follows
from Mineyev-Yaman~\cite{MinYam}, but here we give a direct proof).
Note that any 
compact orientable $n$-manifold $M$ satisfies 
$||M,\d M||\ge ||\d M||/n$ (see e.g.~\cite{Kue-mul}), hence
$||R,\d R||>0$ holds trivially when $||\d R||>0$.

It is easy to see that the fundamental groups of the
manifolds obtained via relative or non-relative
strict hyperbolization split as nontrivial amalgamated products.
In particular, they do not have
Kazhdan property (T), and hence these manifolds are not homotopy
equivalent to a quaternionic or Cayley hyperbolic manifold.
Furthermore, using by now standard harmonic map technology,
we show the following.

\begin{thm} If $K$ is a finite simplicial
complex, then $\H(K)$ is not homotopy equivalent to a 
compact K\"ahler manifold of real dimension $\ge 4$.
\end{thm}

\begin{thm} 
If $R$ is a manifold
obtained by relative strict hyperbolization, then
$\Int (R)$ is not homeomorphic to
an open subset of a compact K\"ahler manifold of 
real dimension $\ge 4$.
\end{thm}

By contrast, it is much harder to decide when $R$ admits
a real hyperbolic metric (see Section~\ref{sec: neg curved}
for some examples).

An important (and still poorly understood) 
invariant of a relatively hyperbolic 
group is the Bowditch boundary, 
which generalizes both the ideal boundary of a hyperbolic 
group, and the limit set of a
geometrically finite isometry group of a negatively pinched
Hadamard manifold. Any relatively hyperbolic group $G$ 
acts on its Bowditch boundary as a geometrically finite
convergence group whose maximal geometrically parabolic subgroups
are exactly the maximal parabolic subgroups of $G$. In fact,
Yaman~\cite{Yaman} characterized relatively hyperbolic
groups as geometrically finite convergence groups, generalizing
a similar characterization of hyperbolic groups due to Bowditch.
For the fundamental groups of complete finite volume negatively
pinched manifolds the Bowditch boundary is a sphere, so
we ask the following.

\begin{quest} Let $R$ be an aspherical manifold of dimension $\ge 4$ 
obtained by the relative strict hyperbolization. 
What is the Bowditch boundary of $\pi_1(R)$? Is it ever a sphere?
Is it a sphere when the maximal parabolic subgroups are virtually nilpotent?
\end{quest}

The structure of the paper is as follows. 
In Section~\ref{sec: strict hyp}-\ref{sec: build block} 
we review the strict hyperbolization,
and then in Section~\ref{sec: strict rel hyper} we
discuss the relative strict hyperbolization. 
Theorem~\ref{intro-thm: rel hyp} is proved in 
Section~\ref{sec: proving rel hyp}. Section~\ref{sec: simpl vol info}
contains background on relative simplicial volume.
Section~\ref{sec: hyp and simp vol} 
discusses relations between the relative strict
hyperbolization and the simplicial volume; 
Theorems~\ref{intro: thm-hausd}--\ref{intro-thm: finiteness}
are proved in this section. 
Section~\ref{sec: cohopf} is devoted to applications to the co-Hopf
property. Theorem~\ref{intro-thm: no split} on non-existence of
elementary splittings is proved in Section~\ref{sec: no split}.
Sections~\ref{sec: kahler}-\ref{sec: neg curved}
discuss when manifolds obtained by relative or non-relative 
strict hyperbolization admit negatively curved Riemannian metrics,
in particular K\"ahler ones.
In Section~\ref{sec: acyl for maps} we prove general acylindricity
results for CW-pairs with relatively hyperbolic fundamental groups.
Section~\ref{sec: geod bound} contains a proof that compact hyperbolic
manifolds with totally geodesic boundary satisfy 
Assumption~\ref{intro: assum}.

\section{Strict hyperbolization}
\label{sec: strict hyp}

The (non-relative) strict hyperbolization is defined 
in~\cite[Theorem 7.6]{CD} as the combination of two independent 
constructions due to Gromov and Charney-Davis.
The Gromov's construction is described 
in~\cite[4c]{DJ}, \cite[Section 15]{Davis-ex}, and~\cite[pp 348--349]{CD}.
It takes as the input an arbitrary $n$-dimensional simplicial complex 
$K$ and turns it into a locally $CAT(0)$ cubical complex denoted by
$\mathcal{G}(K)$. By~\cite[Lemma 7.5]{CD} there is a 
folding map $p\co\G(K)\to\square^n$, where by a {\it folding map}
we mean a cellular map whose restriction on every
cell is a combinatorial isomorphism.

The strict hyperbolization of Charney-Davis~\cite[Proposition 7.1]{CD} 
associates to an $n$-dimensional cubical complex $C$ a piecewise-hyperbolic
$n$-dimensional cell complex $\mathcal S(C)$ such the links of 
the corresponding cells in $C$ and $\mathcal S(C)$ are isometric 
piecewise-spherical polyhedra. Gromov 
showed~\cite[Proposition II.4.14]{BH} that a 
piecewise-Euclidean (or piecewise-hyperbolic) cell complex is $CAT(0)$ 
(or $CAT(-1)$, respectively) if and only if the link of every vertex 
is $CAT(1)$. It follows that if $C$ is $CAT(0)$, then $\mathcal S(C)$ 
is $CAT(-1)$. 

A key idea in constructing $\mathcal S(C)$ is that for each $k$ 
there exists a compact orientable hyperbolic
$k$-manifold with corners $\mathcal X^k$ such that 
the boundary of $\mathcal X^k$ is 
subdivided into totally geodesic $(k-1)$-dimensional 
faces that intersect orthogonally so that the poset 
of faces of $\mathcal X^k$ is the isomorphic to the poset
of faces of $\square^k$.  To obtain $\S(C)$ 
one glues together the $\mathcal X^k$'s in the same combinatorial
pattern as $k$-cubes of $C$. Furthermore, by~\cite[Lemma 5.9]{CD}
for each $k$ there exists a smooth face-preserving degree one 
maps $f_k\co X^k\to\square^k$, and 
if $C$ admits a folding map $p\co C\to\square^n$, one can
describe $\S(C)$ as the subset of $C\times X^n$
consisting of all $(c,x)$ that satisfy $p(c)=f_n(x)$.

We let $\H(K):=\S(\G(K))$, and refer to 
$\H(K)$ as the {\it strict hyperbolization} of $K$. 
The strict hyperbolization has the following properties
which can be deduced from~\cite{DJ, CD, Davis-ex}.
\begin{itemize}
\item
(Negative curvature) $\H(K)$ is locally $CAT(-1)$
piecewise-hyperbolic cell complex. 
If $\dim(K)\le 1$, then $\H(K)$ is the barycentric subdivision of $K$.
\item
(Hyperbolized cell)
If $\s$ is a simplex of dimension $n>0$, then $\H(\s)$ 
is a connected, compact, orientable 
$n$-manifold with boundary. 
\item
(Functoriality) If $i\co J\to K$ is an embedding onto
a subcomplex, then there is an isometric embedding
$\H (i)\co \H(J)\to\H(K)$ onto a locally convex subspace. 
(Charney-Davis~\cite{CD} called locally convex subspaces
``totally geodesic'').
\item 
(Preserving local structure) 
The link of $\H(\s)$ in $\H(K)$ is PL-homeomorphic to the link 
of $\s$ in $K$. Furthermore,
if $v$ is vertex of $K$, then the link of the vertex $h^{-1}(v)$ 
is isomorphic to a subdivision of the 
link of $v$.
\item (Small balls centered at vertices are cones)
If $v$ is a vertex of $\H(K)$, then the link $L_v$ of 
$v$ is a $CAT(1)$ piecewise-spherical complex, and 
by a result of Berestovskii,
a sufficiently small $\e$-ball centered at $v$ is isometric 
to the $\e$-ball centered at the cone point
in the $-1$-cone over $L_v$~\cite[page 59 and Theorem 3.14]{BH}.
Alternatively, the ball can be described as the
warped product $[0,\e)\times_{\sinh(r)} L_v$ 
with base $[0,\e)$, warping
function $\sinh(r)$, and fiber $L_v$. 
(For information on warped products see e.g.~\cite{AB} where
the above warped product is called an {\it elliptic cone} over $L_v$).
\item
(Hyperbolization map)
There is a PL-map $h\co\H(K)\to K$ that restricts to a
degree one map 
$(\H(\s),\partial\H(\s))\to (\s, \partial\s)$ for each simplex $\s$.
One can choose $h$ so that for each vertex $v$ of $K$ it
induces a PL-isomorphism between the small conical neighborhoods
of $h^{-1}(v)$ and $v$. 
\item 
(Surjectivity) $h$ induces surjections on homology and 
fundamental groups.
\item
(Preserving invariants of manifolds) If $K$ is a triangulated 
manifold, or PL manifold, or a smoothly triangulated smooth manifold,
then so is $\H(K)$, and $h$ pulls back the rational Pontrjagin classes,
and the first Stiefel-Whitney class. (This boils down to the fact that 
$\mathcal X^n$ is orientable and has trivial rational Pontrjagin classes). 
In particular, if $K$ is oriented,
then so is $\H(K)$. 
\end{itemize}

\section{On the building block in the strict hyperbolization}
\label{sec: build block}

The manifold $\mathcal X^k$ is obtained from a carefully chosen
closed orientable hyperbolic $k$-manifold by cutting it open 
along a family of $k$ connected totally geodesic $(k-1)$-dimensional 
submanifolds that intersect orthogonally. 
There is a considerably flexibility in choosing $\mathcal X^k$ 
and its structure is not yet well understood. Of course,
$\mathcal X^2$ is just a compact orientable hyperbolic
surface whose boundary is a piecewise geodesic ``square''. 
In higher dimensions
visualizing $\mathcal X^k$ seems much harder, and it is
of some interest to get a better grip on the topology of
$\mathcal X^k$. Here are some open questions:
\begin{itemize}
\item
Can one arrange $\mathcal X^k$ to be stably parallelizable? 
This question was asked in~\cite[p. 348]{CD}, and
if the answer is "yes", then the strict hyperbolization 
construction with this building block preserves 
stable tangent bundle. 
\item
The poset of faces of $\mathcal X^k$ is the poset of faces of 
the $k$-cube, yet for $i>0$ the $i$-dimensional
faces of $\mathcal X^k$ are generally not 
connected. Can one build $\mathcal X^k$ with connected faces?
If yes, visualizing $\mathcal X^k$ would be considerably easier.
\item
Can one arrange $\mathcal X^k$ to be topologically ``complicated'',
e.g. to have all Betti numbers nonzero? This becomes relevant in
Section~\ref{sec: kahler}.
\item 
As we note below $\mathcal X^k$ is never simply-connected.
What is the ``simplest possible'' $\mathcal X^k$?
\end{itemize}

\begin{lem}\label{lem: block nonab free}
If $k\ge 2$, then $\pi_1(\mathcal X^k)$ contains an nonabelian
free subgroup.
\end{lem}
\begin{proof}
By functoriality
of the strict hyperbolization~\cite[p. 347]{CD}, any 
$i$-dimensional face $\square^i\subset\square^k$ 
gives rise to an embedding $\mathcal X^i\subset\mathcal X^k$
onto a locally convex submanifold, where $\mathcal X^i$ is contained in
(but not necessarily equal to) the $i$-dimensional face of $\mathcal X^k$
that corresponds to $\square^i$. Since $\mathcal X^i$ 
a locally convex subspace,
the embedding $\mathcal X^i\subset\mathcal X^k$ is 
$\pi_1$-injective~\cite[Proposition II.4.14]{BH}, and
in particular, $\pi_1(\mathcal X^k)$ contains a nonabelian
free subgroup isomorphic to $\pi_1(\mathcal X^2)$.
\end{proof}

\begin{lem} \label{lem: block embeds hyp mfld}
There exists a cubical complex structure on the
$n$-torus $T^n$ such that $\S(T^n)$ is a closed
hyperbolic manifold, i.e.
a Riemannian manifold of constant negative curvature. In particular,
there is a $\pi_1$-injective embedding of
$\mathcal X^n$ into  closed hyperbolic $n$-manifold.
\end{lem}
\begin{proof}
Let $T^n$ be the $n$-torus built by identifying the opposite sides
of the cube $[-1,1]^n$ in $\mathbb R^n$. Every quadrant of
$\mathbb R^n$ intersects $[-1,1]^n$ in a cube of sidelength $1$,
and these $2^n$ cubes give $T^n$ the structure of a cubical
complex. 
(The definition of a cubical complex that we use 
follows~\cite[p.330]{CD}, and is slightly more general than
the definition in~\cite{BH}: namely, 
we require that each cube in a cubical complex is embedded, 
yet we allow distinct cubes to intersect in several faces, 
rather than in just one face).
Note that the corresponding piecewise-Euclidean metric
on $T^n$ is flat, or equivalently the link at each face is a 
round sphere. Then by~\cite[p.347]{CD} the link at each face of
$\S(T^n)$ is also a round sphere, therefore, the 
piecewise-hyperbolic metric on $\S(T^n)$ is 
in fact hyperbolic. The inclusion $\mathcal X^n$ onto a building block of
$\S(T^n)$ is locally convex, and 
hence $\pi_1$-injective~\cite[Proposition II.4.14]{BH}.
\end{proof}

\begin{lem}\label{lem: retract and T} 
If $K$ is connected and $\dim(K)=n$, then $\H(K)$
retracts onto $\mathcal X^n$.
If $n\ge 2$, then $\pi_1(\H(K))$
does not have Kazhdan property (T).
\end{lem}
\begin{proof}
We think of $\H(K)$ as the
subset of $\G(K)\times \mathcal X^n$
consisting of all $(c,x)$ that satisfy $p(c)=f_n(x)$,
where $p\co\G(K)\to \square^n$ is a folding map, and
$f_n\co \mathcal X^n\to\square^n$ is the map of~\cite[Lemma 5.9]{CD}.
Fix an $n$-cube $i\co\square^n\hookrightarrow\G(K)$ 
in $\G(K)$, which exists as $\dim(K)=n$, and let $\s=i(\square^n)$.
Then 
$\H(K)\cap (\s\times \mathcal X^n)$ is a copy of $\mathcal X^n$ 
that can be thought of as the graph
of $i\circ f\co \mathcal X^n\to \s$. 
Restrict to $\H(K)$ the projection of
$\G(K)\times \mathcal X^n$ onto the second factor, 
and then compose it with the inclusion 
$\mathcal X^n\to \s\times\mathcal X^n$
given by $x\to (x,i(f(x)))$.
This defines a retraction of $\H(K)$
onto $\H(K)\cap (\s\times \mathcal X^n)$, which is a 
copy of $\mathcal X^n$.
Since retractions are $\pi_1$-surjective,
Lemmas~\ref{lem: block nonab free}-\ref{lem: block embeds hyp mfld} 
show that $\pi_1(\H(K))$ surjects
onto the nontrivial group $\pi_1(\mathcal X^n)$
that acts freely on the hyperbolic $n$-space. Thus
$\pi_1(\H(K))$ cannot have property (T)~\cite{dlHarVal}.
\end{proof}

\section{Relative strict hyperbolization}
\label{sec: strict rel hyper}

In this section we adapt the discussion~\cite[Section 2]{DJW} 
of the relative (nonstrict) hyperbolization to the strict case.

Let $K$ be a connected finite-dimensional simplicial complex, 
and let $L$ be a (not necessarily connected, nonempty) subcomplex of $K$.
Let $P$ be the simplicial complex obtained by attaching
a cone $CL_i$ to each path-component $L_i$ of $L$. 
Denote by $o_i$ the cone point of $CL_i$.  
Let $h\co \H(P)\to P$ be the strict hyperbolization map as 
in Section~\ref{sec: strict hyp}.
By functoriality, $\H(K)$ is a subcomplex of $\H(P)$, and
none of the vertices $h^{-1}(o_i)$ lies in $\H(K)$, 
because again by functoriality $h(\H(K))=K$.
Recall that we have chosen $h$ so that each $h^{-1}(o_i)$ 
has some small open $\e$-neighborhood $O_i$ such that  
$h\co O_i\to h(O_i)$ is a homeomorphism. 
Let \[
R_K:=\H(P)\setminus\left(\bigcup_i O_i\right)
\] and let $R_L$ be the (topological) 
boundary of $R_K$ in $\H(P)$, which is also the union of
the boundaries of $O_i$'s.
Since the boundary of $O_i$ is a subdivision of $L_i$, 
we conclude that $R_L$ is a subdivision of $L$.

We refer to the pair $(R_K, R_L)$ as the {\it relative
strict hyperbolization of $(K,L)$}. For future use we
note basic properties of the pair $(R_K, R_L)$.

\bf (A)\rm\ Since $h\co O_i\to h(O_i)$ is a homeomorphism, 
contracting along cone directions
defines a deformation retraction and a homeomorphism
of $(h(R_K), h(R_L))$ onto $(K,L)$.

\bf (B)\rm\  The restriction of $h$ to the map 
$R_K\to h(R_K)$ is surjective on homology 
and fundamental groups, because $h(R_K)$ deformation
retracts onto $K$ and $h\co \mathcal H(K)\to K$ is surjective 
on homology and fundamental groups.

\bf (C)\rm\ The facts C1, C2, C3 below holds for triangulated, 
PL and smooth manifolds (where in the smooth case one works
with smooth triangulations), but we just
state all results in PL category; other cases are similar. 
Let $K$ be a PL manifold with boundary $L$, and let 
$P$ be the simplicial complex obtained by attaching
a cone $CL_i$ to each path-component $L_i$ of $L$.
By~\cite[p 348]{CD} and~\cite[pp 356--357]{DJ}, 
$\mathcal H(P)$ can be identified with $\phi^{-1}(0)$,
where $\phi\co P\times\mathcal X^n\to\R^n$ is the difference 
of the projections $P\to\square^n$, $X^n\to\square^n$,
where $\phi$ is transverse to $0$, and where the restriction
of the projection $P\times\mathcal X^n\to P$ to $\phi^{-1}(0)$ 
corresponds to the map $h\co\mathcal H(P)\to P$.
Restricting $\phi$ to $K\times\mathcal X^n$, we get the map
$\phi\co K\times\mathcal X^n\to\R^n$ where $R_K$ is
identified with $\phi^{-1}(0)$.
Now just like in~\cite[p 348]{CD} and~\cite[pp 356--357]{DJ}
we deduce the following.

\bf (C1)\rm\ 
$R_K$ is a PL manifold with boundary $R_L$, because $\phi$ is 
transverse to $0$.

\bf (C2)\rm\ 
The normal (block) bundle of $R_K$ in $K\times\mathcal X^n\to\R^n$
is trivial, so that the tangent bundle $TR_K$
to $R_K$ is stably isomorphic to the restriction
of $TK\times T\mathcal X^n$ to $R_K$.
Since $\mathcal X^n$ is orientable, and has trivial rational Pontrjagin
classes, the map $h\co R_K\to h(R_K)$ of compact
manifolds pulls back the first Stiefel-Whitney class
and the rational Pontrjagin classes.

\bf (C3)\rm\ 
If the manifold $(K,L)$ is orientable, then so is $(R_K,R_L)$,
because $h$ pulls back the first Stiefel-Whitney class.

\bf (D)\rm\ 
The facts D1, D2, D3 below are 
proved in~\cite[Section 2]{DJW} for nonstrict relative
hyperbolization, and the same proof works verbatim 
if throughout the proof piecewise-Euclidean complexes 
are replaced with piecewise hyperbolic complexes.

\bf (D1)\rm\ The inclusion of each connected component of 
$R_L$ into $R_K$ is $\pi_1$-injective, 

\bf (D2)\rm\ $R_K$ is aspherical if and
only if each component of $R_L$ is aspherical.

\bf (D3)\rm\ 
For a small positive $\delta$,
the $\delta$-neighborhood of $R_L$ in $R_K$ is the warped 
product $[\e,\e+\delta)\times_{\sinh(r)} R_L$, where $R_L$ 
is given the $CAT(1)$ piecewise-spherical metric induced
by the inclusion $R_L\subset\H(K)$.
If $q\co\tilde R_K\to R_K$ is the universal cover
with pullback metric, then each component $C_L$ of
$q^{-1}(R_L)$ is a $CAT(1)$
piecewise-spherical complex, whose $\delta$-neighborhood is
the warped product  $[\e,\e+\delta)\times_{\sinh(r)} C_L$.
Therefore, the space $\bar R_K$ obtained by attaching the
elliptic cone $[0,\e+\delta)\times_{\sinh(r)} C_L$
along $(\e,\e+\delta)\times_{\sinh(r)} C_L$ is locally 
$CAT(-1)$; hence $\bar R_K$ is $CAT(-1)$ because $\bar R_K$
is simply-connected.

\begin{rmk}
Since $h_\ast$ is surjective on homology,
$h^\ast$ is injective on rational cohomology, and hence
if $K$ is a manifold that has a
nontrivial rational Pontrjagin class, $R_K$ has
the same property. Note that by replacing $K$ with the
connected sum of $K$ and a manifold with a nonzero rational
Pontrjagin class, one can always arrange that the
rational Pontrjagin class of $K$ is nonzero. 
Also if $n\ge 4$, then by attaching handles one can prescribe 
the fundamental group of $K$ to be any given finitely presented group,
so one can arrange that $\pi_1(R_K)$ surjects on 
any given finitely presented group. Thus, by varying $K$
while keeping $L$ fixed, one can produce $R_K$ with boundary $L$ and 
fairly complicated topology.
\end{rmk}

\begin{rmk} \label{rmk: glue all ends}
There is another natural way to define the relative 
strict hyperbolization of $(K,L)$, where instead of 
$\mathcal H(K\cup_i CL_i)$ we look at
$\mathcal H(K\cup CL)$, and then remove a small
conical neighborhood of the vertex $\mathcal H (o)$, which
is the preimage of the cone vertex $o$ of $CL$
under the hyperbolization map
$\mathcal H(K\cup CL)\to K\cup CL$. 
As we now explain this new procedure yields the same
pair $(R_K, R_L)$ as above.
Indeed, by functoriality $\mathcal H(K\cup CL)$ is
$\mathcal H(K)$ with $\mathcal H(CL)$ attached along
$\mathcal H(L)$, and in turn, $\mathcal H(CL)$ is 
the union of $\mathcal H(CL_i)$'s
that are all identified at the vertex $\mathcal H(o)$.
Thus $\mathcal H(K\cup CL)$ is obtained from $\mathcal H(K\cup_i CL_i)$ 
by identifying all $h^{-1}(o_i)$'s to the single point 
$\mathcal H(o)$. In other words, $\mathcal H(K\cup CL)$
is obtained from $R_K$ by attaching a cone over $R_L$,
in particular, $\mathcal H(K\cup CL)$ is homeomorphic to
$R_K/R_L$.

Even more generally, the same proof shows the following.
Suppose we partition the set $\{L_i\}$ of path-components of $L$
into the disjoint subsets $S_1,\dots, S_r$, and then attach $r$ 
cones to $K$ with cone points $o_k$, $k=1,\dots ,r$,
so that the cone with vertex $o_k$ is attached over 
the union of $L_i$'s that belong to $S_k$. Then after 
applying $\mathcal H$, we get a locally $CAT(-1)$ complex 
with cone points $h^{-1}(o_1), \dots, h^{-1}(o_r)$,
and after removing the cone points we get $R_K$.
\end{rmk}

\begin{lem}\label{lem: rel hyp, retract, T}
If $\dim(L)<\dim(K)=n$, then $R_K$
retracts onto $\mathcal X^n$.
If $n\ge 2$, then $\pi_1(R_K)$
does not have Kazhdan property (T).
\end{lem}
\begin{proof} 
Note that $\H(K\cup CL)$ has dimension $\le n$ as
$\dim(L)<\dim(K)=n$. 
Since there is an $n$-simplex that does
not lie in $L$, there exists an
inclusion $i\co\mathcal X^n\to\mathcal H(K\cup CL)$ 
such that $i(\mathcal X^n)$ does not contain the
cone vertex $\mathcal H(o)$. 
Hence the retraction of
$\mathcal H(K\cup CL)$ onto $i(\mathcal X^n)$
given by Lemma~\ref{lem: retract and T} restricts
to a retraction of $R_K$ onto $i(\mathcal X^n)$.
Finally, as in Lemma~\ref{lem: retract and T} we
conclude that $\pi_1(R_K)$ does not have property (T).
\end{proof}

\section{Proving relative hyperbolicity}
\label{sec: proving rel hyp}

One of the equivalent definitions of a relatively
hyperbolic group is due to Bowditch~\cite[Definition 2]{Bow-rel}. 
Let $P_1,\dots, P_k$ be finitely generated 
infinite subgroups of a group $G$. 
Let $\mathcal P$ be the set of subgroups of $G$ that are
conjugate to $P_i$, for some $i=1,\dots, k$. Then $G$ is called
{\it hyperbolic relative to} $P_1,\dots, P_k$ if
$G$ admits an action on a connected fine hyperbolic graph $\G$ 
such that each edge has finite stabilizer,
there are only finitely many orbits of edges, and elements
of $\mathcal P$ are precisely the stabilizers of vertices
of infinite valency.
Here a graph is called {\it fine} if for each $L>0$
every edge of $\G$ is contained in only finitely many circuits
of length $\le L$. If in the above situation $P_1,\dots, P_k$ 
are not specified, we shall say, with a slight abuse of language, 
that $G$ {\it hyperbolic relative to} $\mathcal P$. 
Elements of $\mathcal P$
are called {\it maximal parabolic} subgroups; any subgroup of a maximal
parabolic subgroup is called {\it parabolic}.

\begin{ex}
If all vertices of $\G$ have finite
valency (i.e. $\mathcal P$ is empty), then $G$ is quasi-isometric
to $\G$ and hence $G$ is a hyperbolic group. 
\end{ex}

\begin{ex} If $G$ is the fundamental group of a finite graph of groups
whose edge groups are finite, then the Bass-Serre theory gives a
$G$-tree $\G$ with finite edge stabilizers, and $\G/G$ is the original 
graph. Since any tree is hyperbolic and fine, $G$ is hyperbolic relative
to the collection of infinite vertex stabilizers.
\end{ex}

\begin{proof}[Proof of Theorem~\ref{intro-thm: rel hyp}]
Since $R_L$ is a subdivision of a finite complex $L$, 
we know that $R_L$ has finitely many (connected) 
components $L_i$, $i=1,\dots ,r$,
where we assume that $\pi_1(L_i)$ is infinite 
if and only if $i=1,\dots, m$ where $m\le r$.
Each inclusion $L_i\to R_K$ is $\pi_1$-injective
hence it defines a conjugacy class of subgroups of 
$G:=\pi_1(R_K)$ which are isomorphic to $\pi_1(L_i)$, and
let $H_i$ be a subgroup in the conjugacy class.
Let $\mathcal P$ be the union of conjugacy classes
that contain {\it infinite} $H_i$'s.

We need to show that $G$ is hyperbolic relatively to 
$\mathcal P$, and that $G$ is non-elementary.
Look at $G$-action on the universal cover $q\co\tilde R_K\to R_K$.
Conjugates of $H_i$ are in one-to-one correspondence
with components of $q^{-1}(L_i)$, namely,
each subgroup conjugate to $H_i$ is the stabilizer of 
some components of $q^{-1}(L_i)$. These 
components are simply-connected since each inclusion 
$L_i\to R_K$ is $\pi_1$-injective.
As in Section~\ref{sec: strict rel hyper}
for every $i$ we attach an elliptic cone along each 
component of $q^{-1}(L_i)$, denote the result by $\bar R_K$, 
and let $G$ acts isometrically and simplicially on 
$\bar R_K$ in the obvious way so that every conjugate of $H_i$ 
fixes the vertex of the corresponding cone.
Since $\bar R_K$ is $CAT(-1)$, 
$\bar R_K$ is a hyperbolic metric space.

Denote by $\G$ the $1$-skeleton of $\bar R_K$ and check that
the $G$-action on $\G$ satisfies the Bowditch's definition 
of relative hyperbolicity. 
Since $\G$ is quasi-isometric to $\bar R_K$,
the graph $\G$ is hyperbolic. 
Since $K$ and $L$ are finite complexes, 
each $H_i$ is finitely generated, and $\G$ has
only finitely many orbits of edges.
By construction, elements of conjugacy classes of $H_i$'s 
correspond bijectively to cone points.
The cone points are the only vertices that can have nontrivial
stabilizers, so since no edge joins two cone points,
each edge has trivial stabilizer. 

Now we check that $\G$ is fine, which is the only 
part of the proof that is really new.
By~\cite[Proposition 2.1(F5)]{Bow-rel}
it suffices to show that for each vertex $v$, any infinite set of
vertices that are adjacent to $v$ is unbounded in the induced
metric on $\G\setminus\{v\}$.
The only vertices that can have infinitely many adjacent ones are
the cone points stabilized by some subgroup in $\mathcal P$, 
so we can assume that $v$ is one of such cone points,
and denote the corresponding elliptic cone in $\bar R_K$ 
by $C_v$. Note that $C_v$ is convex
because $C_v$ equals to a metric ball $\bar R_K$ centered at $v$
and $\bar R_K$ is $CAT(-1)$~\cite[Proposition II.1.4]{BH}.
Let $S$ be an infinite set of vertices adjacent to $v$, and
arguing by contradiction,
assume that $S$ is bounded in $\G\setminus\{v\}$,
or equivalently that $S$ is bounded in 
$\G\setminus\mathrm{Int}(C_v)$.
Because $\bar R_K$ is $CAT(-1)$
(in fact $CAT(0)$ is enough), the nearest point retraction
$\bar R_K\to C_{v}$ is
 distance-nonincreasing~\cite[Proposition II.2.4]{BH}, 
in particular,
it does not increase the lengths of curves, and hence
any bounded subset of $\G\setminus\mathrm{Int}(C_v)$
projects to a subset of $\partial C_v$
that is bounded in the induced metric on $\partial C_v$.
Since the retraction restricts to the identity on $S$,
we conclude that $S$ is bounded in the induced metric
on $\partial C_v$, and since $S$ is discrete,
it must be finite, which is a contradiction.

Finally, we show that $\pi_1(R_K)$ is non-elementary. 
If $\pi_1(R_K)$ equals to a parabolic subgroup, 
then $\H(P)$ is simply-connected,
and if $\pi_1(R_K)$ is finite or virtually-$\mathbb Z$, then
so is $\pi_1(\H(P))$. (Recall that $\mathcal H(P)$ is $R_K$ with
cone attached over $R_L$).
Take a $2$-simplex $\s^2$ in $P$.  
By functoriality 
$\mathcal H(\s^2)$ is a locally convex subset of $\H(P)$.
Therefore, by~\cite[Proposition II.4.14]{BH}
the inclusion $\mathcal H(\s^2)\to\H(P)$ is $\pi_1$-injective,
so we conclude that $\pi_1(\mathcal H(\s^2))$ is finite or
virtually-$\mathbb Z$.
On the other hand, $\mathcal H(\s^2)$ is a compact surface
with boundary of negative Euler characteristic, 
so its fundamental group is free nonabelian.
\end{proof}

\section{Simplicial volume basics}
\label{sec: simpl vol info}

Here we review the properties of simplicial volume that are
relevant to what we are going to do next.
We refer to~\cite{Grom-vol-bounded-coh, Ben-Pet}
for details.
The complex $C_*(X)$ of (singular real-valued) chains on a topological
space $X$ has a natural $l^1$-norm given by
$||\sum_{i=1}^s r_i\s_i||=\sum_{i=1}^s |r_i|$.
This norm gives rise to a pseudonorm on the homology of $X$, namely
$||h||=\inf\{||z||\co h=[z]\}$, i.e. the pseudonorm of a homology class
$h$ is the infimum of the norms of the cycles representing $h$.
More generally, for a pair of spaces $(X,Y)$ the relative chain complex
$C_*(X,Y)=C_*(X)/C_*(Y)$ has the quotient $l^1$-norm, and the infimum of the 
norms of the cycles representing a relative homology class is the 
a pseudonorm on $H_*(X,Y)$. 

If $M$ is a compact oriented manifold with (possibly empty)
boundary, the {\it relative
simplicial volume} $||M,\partial M||$ of $(M,\partial M)$ is the pseudonorm of 
the fundamental class $[M,\partial M]$. For non-orientable compact manifold
$M$ with the orientation cover $\tilde M$, one defines
$||M,\partial M||=||\tilde M,\partial\tilde M||/2$. 
For example, if the interior of $M$
carries a complete finite volume hyperbolic metric $g$, then 
$||M,\partial M||=\mathrm{Vol}(M,g)/V_n$, where $V_n$ is the 
maximum of the volumes
of ideal regular simplices in the hyperbolic 
$n$-space~\cite{Thu-notes, Fra}.
%
%
\begin{rmk}
It is proved in~\cite{MinYam} that if $M$ satisfies 
Assumption~\ref{intro: assum}, then any class in $H^{k}(M,\d M)$
is bounded for $k\ge 2$. In particular, the $n$-dimensional
cohomology class dual to the fundamental
class $[M,\d M]$ is bounded, which by a standard 
argument~\cite[p. 278]{Ben-Pet}
implies that $||M,\partial M||>0$.
\end{rmk}

\section{Hyperbolization and simplicial volume}
\label{sec: hyp and simp vol}

Let $K$ be a triangulated compact connected
oriented $n$-manifold with boundary $L$. 
To avoid trivialities we assume $n\ge 2$. 
Denote by $s=s(K)$ the number of $n$-simplices in the triangulation. 
Let $R$ be the the compact oriented $n$-manifold with boundary $\d R$
produced via relative strict hyperbolization of $(K,L)$. 
The goal of this section is to show that 
\[ C_1(n)\ge \frac{||R,\d R||}{s}\ge C_2(n)>0,\] 
for some positive constants $C_1(n)$, $C_2(n)$ 
depending only on $n$.

It will be convenient to replace $K$ with its second 
barycentric subdivision; this increases the number of $n$-simplices
by a factor that only depends on $n$. The number of $n$-simplices
of $K\cup CL$ is the sum of the number of $n$-simplices of $K$
and the number of $(n-1)$-simplices of $L$. Since $K$ is a manifold
with boundary $L$, any
$(n-1)$-simplex of $L$ is a face of a 
unique $n$-simplex of $K$. Thus the number of $n$-simplices
of $K\cup CL$ is $\ge s$ and $\le$ than a factor of $s$
that only depends on $n$.

Let $Z=R/\d R$, let $q\co R\to Z$ be the quotient map, 
and let $S=q(\d R)$. 
Alternatively, $Z$ can be thought of as 
the result of attaching a cone to $R$ along $\d R$.
By Remark~\ref{rmk: glue all ends}, $Z=\mathcal H(K\cup CL)$.
Let $i\co S\to Z$ and $j\co Z\to (Z,S)$ be the inclusions.
Then since $S$ is a point and $n\ge 2$, the long exact sequence
of the pair $(Z,S)$ gives the isomorphism 
\[ i_\ast\co H_n(Z)\to H(Z,S).\]
Denote $[Z,S]=q_\ast [R,\partial R]\in H_n(Z, S)$ and
$[Z]=i_\ast^{-1}([Z,S])$.
Note that $Z$ is an oriented pseudomanifold with the fundamental 
class $[Z]$. 
Denote the simplicial norms of $[Z,S]$, $[Z]$ by 
$||Z,S||$, $||Z||$, respectively, 

Denote by $\mathcal H^n(Z)$ the 
$n$-dimensional Hausdorff measure of $Z$. 
As before $\mathcal X^n$ denotes the compact hyperbolic 
manifold with corners constructed in~\cite{CD}
that is used as a building block in the strict hyperbolization.
If $Z$ is built from $N$ copies of $\mathcal X^n$, then clearly
$\mathcal H^n(Z)=N\cdot \mathrm{Vol}(\mathcal X^n)$.
It follows from the definition of strict hyperbolization
that $s\le N\le C(n)s$ where $C(n)$ is a constant that only depends 
on $n$. This estimate combined with 
Proposition~\ref{prop: estimates on simpl vol} below
yields Theorem~\ref{intro: thm-hausd}.

\begin{prop}
\label{prop: estimates on simpl vol}
There exist constants $C_1(n)\ge C_2(n)>0$
depending only on $n$ such that 
$C_1(n)s\ge ||R,\d R||\ge ||Z,S||\ge \frac{||Z||}{n+2}\ge 
\frac{(n-1)!}{\pi (n+2)}\mathcal H^n(Z)\ge C_2(n)s$. 
\end{prop}
\begin{proof}
Triangulate $\mathcal X^n$, 
look at the induced triangulation of $Z$,
and furthermore pass to the second barycentric
subdivision of $Z$, so that $R$ is homeomorphic to $Z$ with
the open star of each cone vertex removed. 
Note that the number of $n$-simplices of the triangulation of
$R$ is at most the number of $n$-simplices in the triangulation of $Z$
that we just constructed, which is $\le s\cdot C_1(n)$
where $C_1(n)$ is some constant depending only on $n$.

The fundamental class
$[R,\d R]$ is represented by a cycle in which each $n$-simplex
of the above triangulation enters only once. This gives
the first inequality $||R,\d R||\le s\cdot C_1(n)$. 
Since the simplicial norm does not increase
under continuous maps, we get the second inequality:
$||R,\d R||\ge ||q_\ast[R,\d R]||=||Z,S||$. 

Let $i_\#$, $j_\#$ be the chain maps induced by $i, j$ in the
short exact sequence of complexes
\[ 
0\to C_\ast(S)\overset{i_\#}{\to} C_\ast (Z)
\overset{j_\#}{\to} C_\ast (Z,S)\to 0.
\]
Since $C_\ast (Z,S)$ carries the quotient norm,
$||Z,S||$ is the infimum of $||c||$ where
$c$ is an $n$-chain in $Z$ and $j_\#(c)$ is a cycle that represents
$[Z,S]$. Note that $c$ need not be a cycle in $C_n(Z)$. 
The strategy of the proof is to 
modify $c$ into a cycle $c^\prime$ such that 
$j_\#(c^\prime)$ still represents $[Z,S]$, and 
$||c^\prime||\le (n+2)||c||$; then we would
get $||Z||\le ||c^\prime||\le (n+2)||c||$
which gives the third inequality after taking
infimum over all $c$'s.

The cycle $c^\prime$ produced after an elementary diagram chase
which we present below.
Since $j_\#$ is a chain map, 
$j_\#\d (c)=\d j_\#(c)=0$, hence
$\d c\in \ker (j_\#)=\mathrm{Im}(i_\#)$, and so  
$\d c=i_\#e$ for some $e\in C_{n-1}(S)$.
Write $e=s\e^{n-1}$ where $s\in\R$,
 and $\e^k\co\Delta^k\to S$ is the unique
singular $k$-simplex with image $S$.
Now $i_\#\d (e)=\d i_\# (e)=\d\d c=0$, so
by injectivity of $i_\#$ we get $\d e=0$.
Since $H_{n-1}(S)=0$, the cycle
$e$ is the boundary of some chain 
$f=t\e^n\in C_n(S)$,
so that $s \e^{n-1}=t \d\e^n$.
Since $\d i_\# (f)=i_\#\d(f)=i_\# (e)=\d c$, 
we deduce that
$c-i_\# f$ is a cycle, which represents 
$[Z]=j_\ast^{-1}([Z,S])$ because $j_\#(c-i_\# f)=j_\#(c)$
represents $[Z,S]$. 

To compute $f$
notice that $\d \e_\a^n=0$ if $n$ is odd, and
$\d \e_\a^n=\e_\a^{n-1}$ if $n$ is even. 
Thus if $n$ is odd, then $\d\co C_n(S)\to C_{n-1}(S)$
is the zero map, so that $0=\d c=i_\# e$ which implies $e=0$, 
and therefore we may take $f=0$.
If $n$ is even, then 
$s \e^{n-1}=e=\d f=t \d\e^{n}=t \e^{n-1}$, i.e.
$s=t$, and therefore in either case
$||e||=||f||$.
Also since $i_\#$ is norm-preserving,
$||i_\# f||=||f||$ and $||\d c||=||i_\# e||=||e||$,
thus all these chains have equal norms. Therefore,
\[ 
||Z||\le ||c-i_\#f||\le 
||c||+||i_\# f||=
||c||+||f||=||c||+||\d c||\le 
(n+2)||c||,
\]
where the last inequality follows from
the estimate $||\d c||\le (n+1)||c||$ coming
from the definition of the boundary homomorphism.
Taking infimum over all $c$'s we get
$||Z||\le (n+2)\cdot ||Z,S||$.

The fourth inequality follows by a result 
of Yamaguchi~\cite{Yamag} who, using 
Thurston's straightening as in~\cite{Grom-vol-bounded-coh}, 
proved the bound when $Z$ is a
compact locally-$CAT(-1)$ geodesically complete
orientable pseudomanifold. 
That $Z$ is geodesically complete follows 
as in~\cite[Proof of Proposition II.5.12]{BH}
because the local $n$-homology group is nontrivial 
at each point of $Z$.

The fifth inequality follows because 
$\mathcal H^n(Z)=N\cdot \mathrm{Vol}(\mathcal X^n)\ge 
s\cdot\mathrm{Vol}(\mathcal X^n)$
by letting $C_2(n)=\mathrm{Vol}(\mathcal X^n)$.
\end{proof}

\begin{rmk}\label{rmk: nonorient estimates on simpl vol}
If $(R,\d R)$ is a {\it nonorientable} compact manifold
obtained via relative strict hyperbolization, we still get
the estimate 
\[ C_1(n)s\ge ||R,\d R||\ge C_2(n)s,\]
for some
$C_1(n)\ge C_2(n)>0$. Indeed, let $\check R$ be the orientable
$2$-fold cover of $R$; this is a manifold with boundary
with a natural triangulation that is pulled backed from $R$.
If $\check Z$ is the space
obtained by attaching a cone to each component of $\d\check R$,
then $\check Z$ is locally-$CAT(-1)$ 
(see the proof of~\cite[Lemma 2.6]{DJW}). Then the proof
of Proposition~\ref{prop: estimates on simpl vol} applies
to $\check R$ and $\check Z$, so we get the above estimates for
$||\check R,\d \check R||$, and hence for $||R,\d R||$, perhaps
for some other $C_1(n)$, $C_2(n)$.
\end{rmk}

\begin{rmk}
Proposition~\ref{prop: estimates on simpl vol} immediately implies
Theorem~\ref{intro-thm: finiteness}
because there are only finitely many ways to glue finitely
many simplices. The same is true in smooth category by using 
smooth triangulations, and smooth relative hyperbolization. 
\end{rmk}

\begin{rmk}
Since the relative simplicial volume is invariant 
under homotopy equivalences of pairs, 
Theorem~\ref{intro-thm: finiteness} gives a finiteness 
result for manifolds pairs obtained by relative strict 
hyperbolization that are in the same homotopy type of pairs. 
Furthermore, if $R$, $R^\prime$ are homotopy equivalent manifolds
obtained by relative strict hyperbolization, and each component
of $\d R$ is aspherical and has intrinsically elementary
fundamental group, then there is homotopy equivalence of pairs
$(R,\d R)\to (R^\prime,\d R^\prime)$ 
(see the proof of Theorem~\ref{thm: finite coker detailed}).
Hence we deduce the following corollary which also holds 
in the smooth category.
\end{rmk}

\begin{cor} Every homotopy type contains  
at most finitely many pairwise PL-non-homeomorphic aspherical
$n$-manifolds $R$ obtained by the relative strict
hyperbolization and such that each component
of $\d R$ has intrinsically elementary
fundamental group.
\end{cor}

\begin{rmk} 
An analog of Proposition~\ref{prop: estimates on simpl vol}
is also true for manifolds obtained via (non-relative) 
strict hyperbolization. To see this
let $K$ be a closed $n$-manifold triangulated with $q$
simplices of dimension $n$. 
Replace $K$ with its second barycentric subdivision;
this changes $q$ by a factor that only depends on $n$. Let
$St_v$ be an open star of a vertex $v$ of $K$, and let
$Lk_v$ be the link of $v$, which we
think of the topological boundary of $s$. 
Then the relative hyperbolization $(R,\d R)$
of $(K\setminus St_v, Lk_v)$
is obtained by removing from $\H(K)$ 
a small conical neighborhood of the vertex corresponding to $v$,
and in particular, $Z=R/\d R$ equals to $\H(K)$. 
By Proposition~\ref{prop: estimates on simpl vol}
we know that $||Z||/s$ is bounded between two positive
constants that only depend on $n$, where $s$ is the number 
of $n$-simplices in the triangulation of $K\setminus St_v$.
It is straightforward to check that $s/q$ is bounded between
two positive constants depending only on $n$, so the desired
assertion follows.
\end{rmk}

\section{Variations of the co-Hopf property}
\label{sec: cohopf}

The following theorem implies Theorem~\ref{thm: finite coker}
because homeomorphic manifolds have equal relative simplicial volumes.

\begin{thm}
\label{thm: finite coker detailed}
If $M_1$, $M_2$ are compact aspherical $n$-manifolds Satisfying 
Assumption~\ref{intro: assum}, and 
$\phi\co\pi_1(M_1)\to\pi_1(M_2)$ is an injective
homomorphism that maps each maximal parabolic subgroup 
to a parabolic subgroup, then $\phi(\pi_1(M_1))$ has 
finite index in $\pi_1(M_2)$.
If in addition $||M_1,\d M_1||\le ||M_2,\d M_2||$, then
$\phi$ is induced by a homotopy
equivalence of pairs $(M_1,\partial M_1)\to (M_2,\d M_2)$,
so that $\phi$ is onto.
\end{thm}
\begin{proof}
Denote the image of $\phi$ by $\hat G$, and look at the 
$\hat G$-action on the universal covering $\tilde M_2$ 
of $M_2$. 
To prove that $\phi(\pi_1(M_1))$ has 
finite index in $\pi_1(M_2)$, we need to 
show that $\hat M_2:=\tilde M_2/\hat G$ is compact.

Let $H_1,\dots , H_k$ are maximal parabolic subgroups
of $\pi_1(M_1)$ which are in one-to-one correspondence with
path-components $Q_1,\dots , Q_k$ of $\partial M_1$.
Since the group $\phi(H_i)$ is parabolic, it stabilizes a path-component 
$\tilde Y_i$ of $\partial\tilde M_2$.
In fact $\phi(H_i)$ acts cocompactly on $\tilde Y_i$,
else the cohomological dimension of $H_i$ would be $<n-1$.
Each closed aspherical manifold $\hat Y_i:=\tilde Y_i/\phi(H_i)$ 
is an incompressible boundary component of $\hat M_2$,
for if $\hat Y_i$ were compressible, so would be its 
projection to $M_2$. A priori $\d \hat M_2$ could have 
components other than $\hat Y_i$'s. To see that this does
not happen, we double $\hat M_2$ along the union of $Y_i$'s.
The result is an aspherical manifold $\hat D_2$, which has empty
boundary if and only if $\d\hat M_2$ is the union of $Y_i$'s.
The fundamental group of $\hat D_2$ is isomorphic to
the fundamental group of the double $DM_1$ of $M_1$ along
$\d M_1$, hence $\hat D_2$  and $DM_1$ are homotopy equivalent.
So by looking at top-dimensional $\Z_2$-homology, we get
$H_n(\hat D_2;\Z_2)\cong H_n(DM_1;\Z_2)\cong\Z_2$, which
implies that $\hat D_2$ is a closed manifold, and in particular $\hat M_2$
is a compact manifold and $\d \hat M_2$ is the union of $\hat Y_i$'s.
This proves that $\phi(\pi_1(M_1))$ has 
finite index in $\pi_1(M_2)$.

Next we show that $\phi$ induces a homotopy 
equivalence of pairs $(M_1, \d M_1)\to (\hat M_2,\d \hat M_2)$.
The proof relies on
the following standard result which can be 
found e.g. in~\cite[Theorem 8.1.9]{Spa}: given connected
CW-complexes $(X,x)$, $(Y,y)$ where $Y$ is aspherical,
there is a natural bijective correspondence between homotopy 
classes of maps $X\to Y$ and the
conjugacy classes of the induced homomorphisms 
$\pi_1(X,x)\to\pi_1(Y,y)$.
 
Since $\hat M_2$ is aspherical, $\phi$ induces 
a homotopy equivalence $f\co M_1\to \hat M_2$,
as well as homotopy equivalences $f_i\co Q_i\to\hat Y_i$, 
where $f_i$ and $f_{|Q_i}$ are homotopic for each 
$i=1, \dots, k$.
One can use the homotopy of $f_i$ and $f_{|Q_i}$
to modify $f$ on a collar neighborhood of $\d M_1$ 
so that we get a map of pairs
$f\co (M_1,\partial M_1)\to (\hat M_2,\partial \hat M_2)$
that induces $\phi$, and equals to $f_i$ when restricted to $Q_i$. 
The same argument applied to a homotopy inverse of $f$ 
yields a map of pairs
$g\co (\hat M_2,\partial \hat M_2)\to (M_1,\partial M_1)$ inducing
$\phi^{-1}$. The composition $g\circ f$ induces the
identity on $\pi_1(M_1)$ and on each $H_i$.
In particular, $g\circ f$
is homotopic to the identity as a selfmap of $M_1$, and each $Q_i$.
Yet we need a stronger conclusion that $g\circ f$ is homotopic to
the identity as a selfmap of $(M_1, \d M_1)$. 
As before, we modify $g\circ f$ on a collar neighborhood
of $\d M_1$ so that $g\circ f$ becomes a selfmap of $(M_1, \d M_1)$
that restricts to the identity of $\d M_1$.
Now by a standard fact~\cite[Proposition 0.19]{Hat}, 
$g\circ f$ is homotopic to the identity 
$rel\ \d M_1$.  
The same argument applies to $f\circ g$. Thus $f$ is a homotopy
equivalence of pairs.

Now we are ready to show that $\phi$ is onto provided 
$||M_2,\d M_2||\ge ||M_1,\d M_1||$. 
There are three cases to consider as follows.

If $M_2$ is orientable, then so is $\hat M_2$. 
Since any homotopy equivalence of 
pairs preserves the Stiefel-Whitney 
classes~\cite[Theorem 6.10.7]{Spa},
$M_1$ is orientable. Furthermore,
homotopy equivalent pairs have equal
relative simplicial volumes:
$||M_1,\d M_1||=||\hat M_2,\d \hat M_2||$.
Finally, since $\hat M_2\to M_2$ 
is a finite cover of degree $|G:\hat G|$, we get
$||\hat M_2,\partial \hat M_2||\ge 
|G:\hat G|\cdot ||M_2,\partial M_2||$.
Combining this information, we get 
$||M_1,\d M_1||\ge |G:\hat G|\cdot ||M_1,\d M_1||$
which forces $|G:\hat G|=1$, because 
$||M_1,\d M_1||>0$ by Theorem~\ref{intro: thm-hausd}.

If $M_2$ and $\hat M_2$ are not orientable, then 
the covering $q\co \hat M_2\to M_2$ does not factor through
the orientation two-fold cover $M_2^\prime\to M_2$, so the $q$-pullback
of $M_2^\prime\to M_2$ is a certain two-fold cover 
over $\hat M_2$, which is orientable, because it is 
a cover of $M_2^\prime$.
Since the orientation two-fold cover of a 
manifold is unique, the two-fold cover over $\hat M_2$ defined
in the previous sentence
has to coincide with the orientation cover $\hat M_2^\prime\to\hat M_2$.
Thus $q$ lifts to a map of the
orientation two-fold covers 
$q^\prime\co \hat M_2^\prime\to M_2^\prime$.
Again by~\cite[Theorem 6.10.7]{Spa}, the
non-orientability of $\hat M_2$ implies that $M_1$ 
is not orientable and that
$f$ lifts to the homotopy equivalence 
$f^\prime\co (M_1^\prime, \d M_1^\prime)\to 
(\hat M_2^\prime ,\d \hat M_2^\prime)$
of the two-fold orientation covers of 
$(M_1,\d M_1)$ and $(\hat M_2, \d\hat M_2)$.
Now the argument of the previous paragraph
applied to $M_1^\prime$, $\hat M_2^\prime$, $M_2^\prime$, 
implies that $q^\prime$ is one-to-one, and hence so is $q$.

Finally, if $M_2$ is not orientable, while $\hat M_2$ is 
orientable, then the covering $q\co \hat M_2\to M_2$ factors 
through the orientation two-fold cover $M_2^\prime\to M_2$.
Since the covering $\hat M_2\to M_2^\prime$ has degree
$|G:\hat G|/2$, we conclude that
\[
||\hat M_2,\partial \hat M_2||\ge 
||M^\prime_2,\partial M^\prime_2||\cdot |G:\hat G|/2=
||M_2,\partial M_2||\cdot |G:\hat G|,
\] where the last equality
holds since $||M^\prime_2,\partial M^\prime_2||=2||M_2,\partial M_2||$.
Now the argument is concluded as in the case when $M_2$ is orientable.
\end{proof}

\section{No splitting over elementary subgroups}
\label{sec: no split}

\begin{proof}[Proof of Theorem~\ref{intro-thm: no split}]
The assertion of (i) was proved 
in~\cite[Section 3]{Bel-Top} in case $M$ is 
negatively pinched, and the same proof works verbatim 
if one replaces ``virtually nilpotent'' by ``elementary'' 
throughout the proof. For completeness we now review the argument
in~\cite[Section 3]{Bel-Top}. 

The proof 
starts off by making the splitting ``relative'', and this is done 
in~\cite[Lemma 3.1]{Bel-Top},
that is essentially due to Bowditch, which applies provided 
every maximal parabolic subgroup of  $\pi_1(M)$ is one-ended.
One-endedness holds in the setting of Theorem~\ref{intro-thm: no split}
because any maximal parabolic subgroup is the fundamental
group of a component of $\d M$, which is a closed aspherical
manifold of dimension $>1$, and any such group is one ended,
as can be seen via a straightforward Mayer-Vietoris argument in 
group cohomology. 
After relativizing the splitting, we
double $M$ along appropriate part of its boundary, 
and then again look at the Mayer-Vietoris
sequence of the splitting. Then a fairly tedious and lengthy
elementary cohomological arguments lead to a contradiction,
which implies (i).

To prove (ii) note that according to~\cite[Theorem 1.14]{DruSap-rips} if 
a non-elementary relatively hyperbolic group $G$ is {\it not} co-Hopfian, 
then either $G$ is isomorphic to a parabolic subgroup, or $G$
splits over an elementary subgroup. The latter possibility
is ruled out by (i). Arguing by contradiction, consider 
a component $B$ of $\d M$ and a monomorphism
$\phi\co\pi_1(M)\to\pi_1(B)$. Since $B$ is a closed aspherical manifold,
$\phi(\pi_1(B))$ has finite index in $\pi_1(B)$. Since
$\phi(\pi_1(M))$ sits between $\phi(\pi_1(B))$ and $\pi_1(B)$,
the index of $\phi(\pi_1(B))$ in $\phi(\pi_1(M))$ is finite,
so that $\pi_1(B)$ has finite index in $\pi_1(M)$. Then the 
intersection of all conjugates of $\pi_1(B)$ is
a finite index subgroup of  $\pi_1(B)$ that is normal
in $\pi_1(M)$. But a maximal parabolic subgroup always
contains the normalizers of its infinite subgroups, so
$\pi_1(M)=\pi_1(B)$, contradicting the assumption that 
$\pi_1(M)$ is non-elementary. 

Finally, we prove (ii). According to~\cite[Corollary 1.14]{DOS}
if $G$ is hyperbolic relative to subgroups
$P_i$ and each $P_i$ is hyperbolic relative to
subgroups $P_{ij}$, then $G$ is hyperbolic relative to $P_{ij}$.
Since the fundamental group of each 
component of $\d M$ has property (m), we can assume
that $\pi_1(M)$ is hyperbolic relative to a family of
subgroups $Q_l$ where each $Q_l$
is contained in the fundamental group of
a component of $\d M$, and furthermore, no $Q_l$
is isomorphic to a non-elementary relatively hyperbolic group. 
By~\cite[Theorem 1.12]{DruSap-rips} if 
a non-elementary relatively hyperbolic group $G$ has infinite 
outer automorphism group and if no maximal parabolic subgroup
of $G$ is isomorphic to a non-elementary relatively 
hyperbolic group, then $G$ splits over an elementary 
subgroup. Thus (iii) follows from (i).
\end{proof}

\begin{rmk} In the case $M$ is negatively pinched, the part (i) of
Theorem~\ref{intro-thm: no split} could be also 
deduced from results of Bowditch that relate the existence of
splittings with absence of global cut points on the boundary,
and this depends on the facts that $n>2$, and that
the Bowditch boundary of $\pi_1(M)$ is the $(n-1)$-sphere in 
the negatively pinched case.
This argument of Bowditch is much more conceptual than the
above proof.
Unfortunately, it does not work in the generality of 
Theorem~\ref{intro-thm: no split}, because in general we know 
nothing about the Bowditch boundary of $\pi_1(M)$. 
\end{rmk}

\section{Strict hyperbolization and K\"ahler manifolds}
\label{sec: kahler}

The following elementary proposition helps to prove that the
manifolds obtained by relative (as well as by non-relative)
strict hyperbolization are not K\"ahler.

\begin{prop} \label{prop: cd1 split}
Let $K$ be a finite connected simplicial complex
of dimension $n\ge 2$ that contains at least one $n$-simplex 
$\s$, and such that not every $2$-simplex of $K$ is a face of $\s$. 
Then $\pi_1(\H(K))$ splits as 
\[
\pi_1(\H(\s))\ast_{\pi_1(\H(\d\s))}\pi_1(\H(K\setminus\Int(\s)))),\]
and $\pi_1(\H(\d\s))$ has infinite index in 
both $\pi_1(\H(\s))$ and $\pi_1(\H(K\setminus\Int(\s)))$.
In particular, $\pi_1(N)$ does not have Serre's property FA.
\end{prop}

\begin{proof} 
The space $\H(K)$ is the union of
$\H(\s)$ and $\H(K\setminus\Int(\s))$ intersecting
along $\H(\d\s)$, which is locally-convex
(as $\d\s$ is a subcomplex of $K$), and hence $\H(\d\s)$
is $\pi_1$-injectively embedded into 
$\H (K)$~\cite[Proposition II.4.14]{BH}. 
Since $\H(K)$ and $\H(\d\s)$ are path-connected, so
is $\H(K\setminus\Int(\s))$.
So $\pi_1(\H(K))$ is an amalgamated product 
of $\pi_1(\H(\s))$ and $\pi_1(\H(K\setminus\Int(\s))))$ 
over $C=\pi_1(\H(\d\s))$. 

Let up show that the index of $C$ in 
$A=\pi_1(\H(\s))$ is infinite.
Since $\H(\s)$ is oriented $n$-manifold with boundary 
$\H(\d\s)$, the homology boundary map is an isomorphism 
in the dimension $n$, so that by exactness $\H(\d\s))\to \H(\s)$ is zero in 
the $(n-1)$-homology. 
If the index $|A:C|$ were finite, then since both $\H(\d\s)$ 
and $\H(\s)$ are aspherical, 
$\H(\d\s)$ would be homotopy equivalent to a finite cover
of $\H(\s)$, so the map $\H(\d\s))\to \H(\s)$ would induce
in the $(n-1)$-homology the multiplication by the degree in the 
of the cover, while by the previous sentence the map is zero.

Let $B=\pi_1(\H(K\setminus\Int(\s)))$, and 
let $\tau$ be a $2$-simplex that
is not a face of $\s$. Then $\pi_1(\H(\d\s\cup\tau))$
is a subgroup of $B$ in which $C=\pi_1(\H(\d\s))$ has infinite
index because each component of 
$\H(\d\s\cap\tau)=\d\s\cap\tau$ is a simplex
of dimension $\le 1$, and $\pi_1(\H(\tau))$ is infinite. Thus
the index $|B:C|$ is infinite.
\end{proof}

\begin{rmk}
A similar argument shows that if 
$K$ is a finite connected simplicial complex of dimension $n\ge 2$,
and $L$ is a subcomplex satisfying $\dim(L)<\dim(K)$, 
then for any $n$-simplex $\s$ of $K$, the group
$\pi_1(R_K)$ splits as an amalgamated product over
over $\H(\d\s)$, where $\H(\d\s)$ has infinite index in each 
of the factors.
\end{rmk}

By a {\it compact K\"ahler manifold}
we mean a closed Riemannian manifold with a K\"ahler metric.

\begin{thm} \label{thm: comp kahler} If $K$ is a finite simplicial
complex, then 
$\H(K)$ is not homotopy equivalent to a compact K\"ahler manifold
of real dimension $\ge 4$.
\end{thm}
\begin{proof}
Arguing by contradiction assume that $\H(K)$ is
homotopy equivalent to a
compact K\"ahler manifold $M$ of real dimension $n\ge 4$. 
First note that 
$\dim(K)=n$ and that $K$ has at least two $n$-simplices.
(Indeed, $\dim(K)=dim(\H(K))\ge n$ because $\H(K)$ is homotopy equivalent
to a closed manifold $M$. If $\dim(K)>n$, then $K$ has a simplex $\s$ 
of dimension $n+1$, 
so $\pi_1(M)\cong\pi_1(\H(K))$ must contain the fundamental groups 
of the closed aspherical $n$-manifold $\pi_1(\d\s)$.
Since $\H(K)$ is aspherical, so is $M$, therefore
$\pi_1(\d\s)$ must have finite index in $\pi_1(M)$ which is impossible
as $\pi_1(\d\s)$ has infinite index in $\pi_1(\s)\le\pi_1(\H(K))$
by the proof of Proposition~\ref{prop: cd1 split}.
If $K$ has only one $n$-simplex $\s$, then $\H(K)$ is the union of 
$\H(\s)$ and a subcomplex of dimension $<n$. Then
$\H(\s)$ and $(\H(K),\H(\s))$ have trivial $n$th homology,
and hence so does $\H(K)$ by the long exact sequence of
the pair which contradicts to the assumption that $\H(K)$ is 
homotopy equivalent to a closed $n$-manifold). 

Hence by Proposition~\ref{prop: cd1 split},
$G=\pi_1(\H(K))$ splits over $\pi_1(\H(\d\s))$.
The splitting $G=A\ast_C B$ gives rise to a nontrivial 
$G$-action on the corresponding Bass-Serre tree $T$, which
has no fixed end. (If $G$ fixes an end, then $A$ fixes the end
and a vertex $v$ of $T$, so that $A$ fixes the edge adjacent to 
$v$ that lies on a unique ray joining $v$ and the end. Thus
the inclusion $C\to A$ is onto contradicting nontriviality of the 
splitting).
Therefore,
by~\cite[Corollary 2.3.2]{KorSch} there exists a $G$-equivariant
harmonic map $u\co\tilde M\to T$, where $\tilde M$ is the universal cover
of $M$. In fact, the image of $u$ locally lies in a finite subtree
of $T$~\cite[p.243]{GroSch} and~\cite[Theorem 3.9]{Sun}, so 
$u$ is pluriharmonic by~\cite[Theorem 7.3]{GroSch}.
Then by the factorization theorem in~\cite[Section 9]{GroSch},
$u$ factors through the hyperbolic plane $\bf H^2$ equipped 
with an isometric discrete (possibly ineffective) 
$G$-action, and moreover the map $\tilde M\to\bf{H^2}$ 
is $G^\prime$-equivariant, where $G^\prime$ is a 
subgroup of $G$ of index $\le 2$.
(By~\cite[Lemma 9.4]{GroSch} 
the discreteness of the $G$-action on $\bf{H^2}$ depends 
on the fact that $G$ is not virtually solvable, which is true 
since $G$ contains free nonabelian subgroup).

The kernel $K^\prime$ of the $G^\prime$-action on $\bf H^2$
lies in the kernel of the $G^\prime$-action on $T$,
which in turn lies in $C$.
The subgroup $C$ is a quasiconvex
subgroup because $\H(\d\s)$ is locally convex in $\H(K)$.
Note that $K^\prime$ is infinite, else $G\cong \pi_1(M)$ 
would have cohomological dimension $\le 2$.
Thus the quasiconvex subgroup $C$ contains $K^\prime$,
which is an infinite normal subgroup of the hyperbolic 
group $G^\prime$, which is impossible by a standard 
``limit set'' argument.
(Indeed, $K^\prime$, $G^\prime$, and $G$ have the same limit
sets, and the limit set of $K^\prime$ lies in the limit
set of $C$. Thus $C$ and $G$ have the same limit set,
hence the quasiconvexity of $C$ implies that $C$ has finite index 
in $G$~\cite{Swe}).
\end{proof}

\begin{rmk} The proof of 
Theorem~\ref{thm: comp kahler} was inspired by
recent work of Delzant-Gromov~\cite{DelGro}
where more general situation is considered
and where the proofs are unfortunately sketchy. 
\end{rmk}

\begin{thm} 
\label{thm: rel hyp Kahler} If $R$ is an $n$-manifold
obtained by relative strict hyperbolization, then
$\Int (R)$ is not homeomorphic to
an open subset of a compact K\"ahler manifold of 
real dimension $\ge 4$. 
\end{thm}

\begin{proof} We think of $R$ as obtained from the 
$\H(K\cup CL)$ by removing a
small conical neighborhood of the cone point. Since
$K$ is an $n$-manifold with boundary $L$, there is 
at least one $n$-simplex $\s$ of $K$ that does not lie in
$L$, and $K\cup CL$ contains a $2$-simplex
that contains the cone point and is not a face of $\s$. 
Thus we get a splitting of 
$\pi_1(\H(K\cup CL))$ as in Proposition~\ref{prop: cd1 split}.

Arguing by contradiction assume that 
$U$ is an open subset of a compact K\"ahler manifold $M$,
and $f\co U\to \Int (R)$ is a homeomorphism.
Note that $f$
extends to a continuous map $\bar f\co M\to R/\d R=\H(K\cup CL)$ that maps
$M\setminus U$ to the cone point $\d R/\d R$, 
because $f$ is proper and continuous, and $R/\d R$
is the one-point-compactification of $\Int (R)$.
The map $\bar f\co M\to R/\d R$ has degree one 
(recall that $R/\d R$ is a pseudomanifold), and 
in particular, $\bar f$ is $\pi_1$-surjective. 

We can write $M$ as the union of $f^{-1}(\H(\s))$ and 
$M\setminus f^{-1}(\Int(\H(\s)))$ which intersect along
$f^{-1}(\H(\d\s))$. Note that $f^{-1}(\H(\d\s))$ is $\pi_1$-injectively 
embedded in both $f^{-1}(\H(\s))$ and $M\setminus f^{-1}(\Int(\H(\s)))$:
indeed, if a homotopically 
nontrivial loop in $f^{-1}(\H(\d\s))$ bounds a disk in $M$,
its $f$-image is a homotopically nontrivial loop in
$\H(\d\s)$ that bounds a disk in $\H(K)$, contradicting the fact that 
$\H(\d\s)$ is $\pi_1$-injectively embedded in $\H(K)$.
Thus by Van Kampen theorem we get a splitting $\bar A\ast_{\bar C} \bar B$
induced by the decomposition of $M$ into $f^{-1}(\H(\s))$ and 
$M\setminus f^{-1}(\Int(\H(\s)))$ with $\bar C$ is the fundamental group
of $f^{-1}(\H(\d\s))$. 
The splitting is nontrivial because it is mapped via $\bar f$
onto the splitting of $R/d R$ along $\H(\d\s)$ that is nontrivial by
Proposition~\ref{prop: cd1 split}.

As in the proof of Theorem~\ref{thm: comp kahler}
we see that $\bar C$ contains an infinite subgroup that is normal
in a subgroup $G^\prime\le \pi_1(M)$ of index $\le 2$. 
Because $f_\ast$ maps $\bar C$ isomorphically onto $\pi_1(\H(\d\s))$,
we conclude that $\pi_1(\H(\d\s))$ contains an infinite subgroup
that is normal in $f_\ast(G^\prime)\le\pi_1(R/\d R)$ 
that has index $\le 2$, because $f_\ast$ is surjective.
The proof is now concluded as in Theorem~\ref{thm: comp kahler}
by looking at limits sets.
\end{proof}

\begin{rmk}
In the statement of
Theorem~\ref{thm: rel hyp Kahler} it should be possible to replace
``not homeomorphic'' by ``not properly homotopy equivalent''. 
This is done in Theorem~\ref{thm: rel hyp Kahler Betti}
under the assumption that
the building block  $\mathcal X^n$ of
the strict hyperbolization has nonzero  
$i$th homology group for some $i>2$. 
\end{rmk}

\begin{ex}
By Theorem~\ref{thm: rel hyp Kahler},
$\Int (R_K)$ is not homeomorphic to a 
quasiprojective manifold (i.e. the complement $M\setminus S$
where $M$ is compact K\"ahler and $S$ is a subvariety with normal 
crossings, see discussion of quasiprojective manifolds 
in~\cite[p.292]{Jos-Yau}). It was shown in~\cite{Yeu}
that any complete finite volume K\"ahler manifold
with negative Ricci curvature and two sided bounds on sectional curvature
is quasiprojective. In particular, any 
complete finite volume K\"ahler manifold $V$ of pinched negative curvature
is quasiprojective, and hence is not homeomorphic to 
$\Int (R_K)$.
\end{ex}

There is somewhat different way to see that $\H(K)$ and 
$R_K$ are not K\"ahler that is also based on harmonic maps
technology, and depends on the assumption that
the building block  $\mathcal X^n$ of
the strict hyperbolization has nonzero  
$i$th homology group for some $i\ge 3$. 
This assumption might be always true if $n\ge 4$
(which is an obvious necessary condition), but I do not
know how to prove it. It is almost certain that
one should be able to {\it choose} $\mathcal X^n$ with nonzero  
$i$th homology group for some $i\ge 3$ (e.g. it would suffice to
find $\mathcal X^n$ that contains two closed totally geodesic
submanifolds of dimensions $i$, $n-i$ that intersect transversely
at one point).

\begin{thm} \label{thm: non-rel hyp Kahler Betti}
Suppose that $\mathcal X^n$ of has nonzero  
$i$th homology group for some $i\ge 3$.
If $K$ is a finite simplicial complex of dimension $n$, 
then $\H(K)$ is not homotopy equivalent to
a compact K\"ahler manifold.
\end{thm}
\begin{proof}
Arguing by contradiction consider
a compact K\"ahler manifold $M$ and a 
homotopy equivalence $M\to\H(K)$.
By Lemmas~\ref{lem: block embeds hyp mfld}-\ref{lem: retract and T}, 
$\H(K)$ retracts onto $\mathcal X^n$, and
$\mathcal X^n$ embeds into a closed hyperbolic
$n$-manifold $Q=\mathcal {S}(T^n)$. Composing the homotopy equivalence
$M\to\H(K)$, the retraction $\H(K)\to\mathcal X^n$, and the
inclusion $\mathcal X^n\to Q$,
we get a map $r\co M\to Q$. By the factorization theorem
of Sampson and Carlson-Toledo~\cite[Theorem 6.21]{ABCKT},
$r$ is homotopic to the map $g\co M\to Q$ that 
factors through a circle or a compact $2$-manifold.
The inclusion $\mathcal X^n\to Q$ is a $\pi_1$-injective
map of aspherical spaces, hence $\mathcal X^n\to Q$ 
lifts to a cover $\bar Q$ of 
$Q$ such that the inclusion $\mathcal X^n\to \bar Q$ is a
homotopy equivalence.
The image of $r$ is $\mathcal X^n$, 
so $r$ lifts to $M\to \mathcal X^n\subset\bar Q$,
which is a homotopy retraction, and therefore is
surjective in the $i$th homology. By assumption
the $i$th homology group of $\mathcal X^n$ is nontrivial,
so $M\to \mathcal X^n\subset\bar Q$ is nontrivial on
the $i$th homology. On the other hand, by the covering
homotopy theorem the map $M\to \mathcal X^n\subset\bar Q$
is homotopic to a lift $\bar g\co M\to\bar Q$
of $g$ that factors through a circle or a compact $2$-manifold,
and is therefore trivial on the $i$th homology as $i\ge 3$.
\end{proof}

\begin{thm} 
\label{thm: rel hyp Kahler Betti} 
Suppose that $\mathcal X^n$ of has nonzero  
$i$th homology group for some $i\ge 3$.
If $L$ is a subcomplex of a finite connected 
simplicial complex $K$ with
$\dim(L)<\dim(K)=n$, then
$\mathrm{Int} (R_K)$ is not properly homotopy equivalent to
an open subset of a compact K\"ahler manifold. 
\end{thm}

\begin{proof}
Arguing by contradiction assume that 
$U$ is an open subset of a compact K\"ahler manifold $M$,
and $f\co U\to \Int (R_K)$ is a proper homotopy equivalence.
Note that $f$
extends to a continuous map $\bar f\co M\to R_K/R_L$ with
$\bar f(M\setminus U)=R_L/R_L$, 
because $f$ is proper and continuous, and $R_K/R_L$
is the one-point-compactification of $\Int (R_K)$.
By Lemma~\ref{lem: rel hyp, retract, T}
there is a retraction $r\co R_K/R_L\to \mathcal X^n$ 
where $\mathcal X^n$ is a block that 
does not contain the cone vertex $R_L/R_L$. 
The composition $r\circ \bar f\co M\to \mathcal X^n$
is surjective on homology, because its restriction 
to $U$ is $r\circ f$ which is surjective on homology as
the composition of a homotopy
equivalence $U\to\Int(R_K)$ and the 
retraction $\Int(R_K)\to X^n$.
The proof is then finished as in the proof of 
Theorem~\ref{thm: non-rel hyp Kahler Betti}.
\end{proof}

\section{Strict hyperbolization and negatively curved manifolds}
\label{sec: neg curved}

In this section we discuss several examples of aspherical manifolds
of the form $\H(K)$ and $\Int(R_K)$ that admit no complete
Riemannian metrics of pinched negative curvature and finite volume.

\begin{ex}\bf (Davis--Januszkiewicz)\rm\
The first example of a closed smoothable
manifold of the form $\H(K)$ that carries no Riemannian 
metric of negative curvature was constructed
in~\cite[p.384]{DJ}. The construction works in each dimension 
$n\ge 5$, and $K$ is the double suspension of 
a smooth homology $(n-2)$-sphere with finite nontrivial 
fundamental group. By Edward's theorem, $K$ is a topological
manifold, yet its obvious triangulation is not PL. 
It s proved in~\cite[p.385]{DJ} that the universal cover
of $\H(K)$ is homeomorphic to $\mathbb R^n$ yet its ideal
boundary as a $CAT(-1)$ space (or equivalently
the ideal boundary of the hyperbolic group $\pi_1(\H(K))$),
is not a sphere, hence $\H(K)$ carries no Riemannian 
metric of negative curvature. (Strictly speaking,
the proof in~\cite{DJ} was written for the Gromov's hyperbolization
procedure that in fact does not yield negative curvature, which was not
known at the time when~\cite{DJ} was written, but
the same proof holds for the strict hyperbolization of~\cite{CD}).
\end{ex}

\begin{rmk}
This kind of examples are impossible in
dimensions $2$ and $3$.
\begin{itemize}
\item Indeed, in dimension $2$, 
the manifolds $\H(K)$ and $R_K$ retract onto $\mathcal X^2$
that has free fundamental group, hence $\H(K)$ and $R_K$
are surfaces of negative Euler characteristics, so their interiors
admit hyperbolic metrics of finite volume. 
\item
In dimension $3$ the manifold $\H(K)$
is hyperbolizable by the Thurston's hyperbolization theorem:
indeed, $\H(K)$ is Haken as it contains the incompressible surface
$\H(\d\s)$ where $\s$ is a $3$-simplex of $K$; also
$\H(K)$ contains no fake $3$-cells because its universal cover is 
$\mathbb R^3$~\cite[Theorem 3b.2]{DJ}, and finally, $\H(K)$ is atoroidal
since its fundamental group is hyperbolic.
\item
Assuming the Poincar\'e Conjecture, one can show that
the aspherical (let us say orientable) 
$3$-manifold $R_K$ is also hyperbolizable, in fact 
this holds for any $3$-dimensional orientable manifold $M$ satisfying 
Assumption~\ref{intro: assum}. Indeed, 
by Theorem~\ref{intro-thm: no split}, the group $\pi_1(M)$ is 
freely indecomposable, so $M$ is irreducible. 
Since $M$ has nonempty boundary, $M$ is Haken. 
Since any $\Z\times\Z$ subgroup of a relatively hyperbolic group
must be parabolic, $M$ is atoroidal.
By Section~\ref{sec: acyl for maps}, $M$ is acylindrical.
Since $\pi_1(M)$ is non-elementary relatively hyperbolic,
it is not virtually abelian, which
excludes various special cases e.g. $M$ is not 
an $I$-bundle over the torus or the
Klein bottle, or that $M$ is not contractible. 
Hence by the Thurston's hyperbolization 
theorem~\cite[Theorem 1.43]{Kap-book} $M$
is homeomorphic to the convex core of a geometrically
finite Kleinian group whose maximal parabolic subgroups
correspond to the boundary tori of $M$. Furthermore,
if $\d M$ consists of tori, then $M$ has the complete
hyperbolic metric of finite volume.
\end{itemize}
\end{rmk}

\begin{ex} 
Long-Reid noted in~\cite{LonRei-eta} that by
the Atiyah-Patodi-Singer
formula the $\eta$-invariant of the cusp cross-section of 
any $1$-cusped hyperbolic $4$-manifold must be an integer,
while there exists a closed orientable flat $3$-manifold 
that has a nonintegral $\eta$-invariant.
However, since any
flat manifold bounds~\cite{HamRoy}, $F$ bounds some $R$
obtained via the relative strict hyperbolization, 
which by~\cite{LonRei-eta} admits no hyperbolic metric.
For related work in complex and quaternion hyperbolic setting
see~\cite{McR}.
\end{ex}

\begin{ex}\label{ex: GT large pinch}
\bf (Gromov-Thurston examples with large pinching)\rm\
Consider an $\Z_n$-action on 
$S^{k+2}:=\mathbb R^{k+2}\cup\{\infty\}$ generated by
a rotation by $2\pi/m$ about $\mathbb R^k\cup\{\infty\}$. 
Fix a $\Z_m$-invariant triangulation of $\overline{\mathbb R^{k+2}}$, 
(and pass to the first barycentric subdivision to ensure 
that the complex can be folded onto a $(k+2)$-simplex). The
result, denoted by $K_m$, is a simplicial complex homeomorphic to the 
$(k+2)$-sphere and equipped with a simplicial semifree
$\Z_m$-action whose fixed point set is homeomorphic to 
the $k$-sphere, which we denote by $S$.
Let $P$  be the subcomplex of $K_m$ 
corresponding to a closed $(k+1)$-dimensional
half-plane in $\mathbb R^{k+2}\cup\{\infty\}$ whose boundary is 
$\mathbb R^k\cup\{\infty\}$. 
An obvious fundamental domain $F_m$ for the $\Z_m$-action
is the ``lens'' bounded by $P$ and 
its image under the generator of $\Z_m$.
Since the strict hyperbolization is functorial, we get 
the induced $\Z_m$-action on $\H(K_m)$ with the fundamental 
domain $\H(F_m)$ and the fixed-point-set $\H(S)$.
Note that $\H(P)$ is a compact manifold whose boundary is $\H(S)$,
and the $\Z_m$-images of $\H(P)$ all intersect
in $\H(S)$. Recall that the strict hyperbolization takes
subcomplexes to subcomplexes that are $\pi_1$-injectively embedded.
Thus $\H(K_m)$ contains the amalgamated product of the fundamental 
groups of the $\Z_m$-images of $\H(P)$ which are amalgamated over
$\pi_1(\H(S))$. If $k\ge 2$, then $\pi_1(\H(S))$
is non-elementary hyperbolic, so the fundamental class of
$\H(S)$ is noncuspidal~\cite[p.10]{GroThu}, yet
it bounds the relative chain given by $P$, 
so we can apply~\cite[p.11]{GroThu} to deduce the following:
for each $n$ there exists a sequence $\{a_m\}$ with
$a_m\to\infty$ as $m\to\infty$ 
such that $\H(K_m)$ is not homotopy equivalent to
a complete Riemannian $n$-manifold with $-a_i^2\le \sec\le -1$. 
Note that we necessarily have $n\ge k+2\ge 4$, because
$\H(K_m)$ is a closed aspherical manifold and $k\ge 2$.
\end{ex}

\begin{ex} 
\label{ex: gt loc symm}
\bf (Gromov-Thurston examples with no locally 
symmetric metric)\rm\
The results in~\cite[p.8--11]{GroThu}, used in 
Example~\ref{ex: GT large pinch},
depend on some delicate volume estimates, and it is 
much easier to see (again following~\cite[p.1--2]{GroThu}) that
$\H(K_m)$ admits no locally symmetric negatively 
curved metric for all large $m$. This is what we do in this
example: the argument is essentially taken 
from~\cite{GroThu} (with a few details added).

Let $h$ be a generator of the group $\Z_m$ 
that acts by simplicial homeomorphisms on $K_m$ and $\H(K_m)$.
Note that $h$ is not homotopic to the identity. (Indeed,
$h$ moves some top-dimensional simplex $\s$ of $K_m$ with 
a vertex $v\in S$ to a different simplex $\s^\prime$, 
so $h(\H(\s))=\H(\s^\prime)$ and $h(v)=v$.
If $\a$ is a homotopically nontrivial loop in $\H(\s)$
based in $v$, then the loops $\a$ and $h(\a)$ are not 
homotopic.) Clearly $h$ is orientation-preserving on $K_m$,
and hence also on $\H(K_m)$.
Suppose $\H(K_m)$ is homotopy equivalent to 
a closed locally symmetric negatively curved manifold
$X/\Gamma$, where $X$ is the corresponding negatively curved
symmetric space and $\Gamma$ is a discrete isometry group
isomorphic to $\pi_1(\H(K_m),\star)$. Let 
$\phi$ be the automorphism of $\Gamma$ induced by $h$.
By Mostow Rigidity, there exists an isometry $\tilde i$
of $X$ such that $\phi(\g)=\tilde i\circ\g\circ\tilde i^{-1}$ for all
$\g\in\pi_1(\H(K_m),\star)$, and $\tilde i$ descends to an isometry
$i$ of $X/\Gamma$ that becomes homotopic to $h$ when 
$\H(K_m)$ is identified with $X/\Gamma$.
Since $h_\ast$ is the identity on $\pi_1(\H(S),\star)$,
we conclude that $\tilde i$ lies in the centralizer
of the subgroup $\Gamma_S$ that corresponds
$\pi_1(\H(S),\star)$. If $\dim(S)\ge 2$, then $\Gamma_S$
is non-elementary, hence $\tilde i$ is a rotation about
the totally geodesic subspace $X_S$ of $X$ that is the convex hull
of the limit set of $\G_S$. 

(In fact, it is easy to see
that $\dim(X_S)=k$ even though we do not need this fact here.
Indeed, $X_S/\Gamma_S$ is an 
aspherical manifold homotopy equivalent to a closed
aspherical manifold $\H(S)$, we get that $\dim(X_S)\ge\dim S=k$.
By the same reason $\dim(X)=\dim(X/\Gamma)=\dim(\H(K_M))=k+2$. 
If $\dim(X_S)=k+2$, then
$\tilde i$ is the identity, contradicting the fact $h$ is not homotopic
to the identity, and if $\dim(X_S)=k+1$, then $\tilde i$
is a reflection in $X_S$, so that $i$ reverses the orientation,
contradicting the fact $h$ is orientation-preserving.)

Since $h^m$ is the identity, $i^m$ is homotopic to the identity,
and hence by Mostow Rigidity, $i^m$ is the identity. We conclude
that $\tilde i$ is an order $m$ rotation about $X_S$.
Therefore, the orbifolds obtained as $\Z_m$-quotients of $X/\Gamma$
are pairwise non-isometric. The volume of any such orbifold is 
$\vol(X/\Gamma)/m$, and $\vol(X/\Gamma)=C(k)||X/\Gamma||$ where $C(k)$ is a constant
depending only on $k$~\cite{Ben-Pet}. 
On the other hand, simplicial volume is invariant under 
homotopy equivalences so $||X/\Gamma||=||\H(K_m)||$ which is bounded
above by (a constant multiple of)
the number of simplices of $K_m$, which in turn
is bounded above by a constant multiple of $m$. Thus there is an 
upper bound depending only on the dimension on the volume of 
the orbifolds obtained as $\Z_m$-quotients of $X/G$. 
By Wang's finiteness theorem there can be only finitely many
isometry types of 
such orbifolds, which shows that for all large $m$, the manifold
$\H(K_m)$ is not homotopy equivalent
to a closed locally symmetric manifold of negative curvature.
\end{ex}

\begin{rmk}
Similarly to Example~\ref{ex: gt loc symm} one can use
relative strict hyperbolization to build examples
of manifolds $R$ such that each component of $\d R$ is 
a flat manifold and $R$ is not homotopy equivalent
to a complete loaclly symmetric negatively curved manifold
of finite volume. 
\end{rmk}

\section{Spaces of maps and acylindricity}
\label{sec: acyl for maps}

It is well-known that two maps into an aspherical space
are homotopic if and only if the induced $\pi_1$-homomorphisms
are conjugate. In this section we take a closer look at
the case when the aspherical space has relatively hyperbolic
fundamental group, and the $\pi_1$-homomorphisms have parabolic
images.

A sample application is that if $M$ be a compact aspherical 
$3$-manifold satisfying Assumption~\ref{intro: assum}, then
$M$ is {\it acylindrical}.

In fact, the argument used to prove acylindricity
has nothing to do with dimension three,
and more generally, instead of $(M,\d M)$ we consider an arbitrary 
CW-pair $(X,Y)$ 
such that \newline
$\bullet$ $X$ is aspherical and locally compact, \newline
$\bullet$ each path-component of $Y$ is aspherical and 
incompressible in $X$, \newline
$\bullet$ the group $\pi_1(X)$ is non-elementary relatively 
hyperbolic, relative to  fundamental groups of 
path-components of $Y$.

Given a connected CW-complex $Z$, and
a continuous map $f\co Z\to Y\subset X$,
we consider the inclusion $\i\co C(Z,Y;f)\to C(Z,X;f)$, 
where $C(Z,Y;f)$ 
denotes the path-component of $f$ in the space of
continuous maps from $Z$ to $Y$ 
with compact-open topology, and similarly for $C(Z,X;f)$.

\begin{prop}\label{prop: acyl}
If $f\co Z\to Y$ is not null-homotopic, then 
the inclusion $\i\co C(Z,Y;f)\to C(Z,X;f)$
is a weak homotopy equivalence. 
\end{prop}

Recall that $\i$ is a weak homotopy equivalence if and only if
for any finite-dimensional CW-pair $(K,J)$ every map
$(K,J)\to (C(Z,X;f),C(Z,Y;f))$ is homotopic $rel\ J$
to a map with image in $C(Z,Y;f)$. For example, applying 
Proposition~\ref{prop: acyl} for $K=[0,1]$, $J=\{0,1\}$,
we conclude that two maps $Z\to Y$ are homotopic if and only if
they are homotopic in $X$.  
In particular if $X$ is a compact 
$3$-manifold with boundary $Y$, and $Z=S^1$, then
$X$ is acylindrical.

\begin{rmk} If $B$ is the component of $Y$ that contains $f(Z)$,
then of course the inclusion $C(Z,B;f)\to C(Z,Y;f)$ is a bijection,
so that $C(Z,X;f)$ is actually weakly homotopy equivalent to  $C(Z,B;f)$,
unless $f$ is null-homotopic. The weak homotopy type of $C(Z,B;f)$
was computed by Gottlieb (see~\cite{Han}), namely
$C(Z, B, h)$ is aspherical (i.e. its universal cover is weakly 
homotopy equivalent to a point), and its
fundamental group is isomorphic to the 
centralizer of $f_\ast\pi_1(Z)$ in $\pi_1(B)$.
Since $f_\ast\pi_1(Z)$ is parabolic, the centralizer of
$f_\ast\pi_1(Z)$ is a parabolic subgroup of $\pi_1(B)$
that contains the center of $\pi_1(B)$.
\end{rmk}

\begin{proof}[Proof of Proposition~\ref{prop: acyl}]
Denote by $\tilde Y$ the preimage of $Y$ under the universal
covering $p\co \tilde X\to X$.
Since components of $Y$ are aspherical and incompressible, 
the components of $\tilde Y$ are contractible.
The identification of $\pi_1(X)$ and the automorphism group 
of the covering $p$ can be chosen to induce a one-to-one correspondence
between the stabilizers in $\pi_1(X)$ of components of $\tilde Y$
(or equivalently, the maximal parabolic subgroups of $\pi_1(X)$),
and the subgroups of $\pi_1(X)$ conjugate to the fundamental groups
of the components of $Y$.
 
Fix an arbitrary path-component $\tilde B$ of $p^{-1}(B)$, and
denote by $P$ the stabilizer of $\tilde B$ in $\pi_1(X)$.
Since $\pi_1(X)$ is torsion-free and relatively hyperbolic,
any two maximal distinct 
parabolic subgroups have trivial intersection. 
In particular,
no element of $P$ can stabilize two distinct
components of $\tilde Y$.
Hence the only non-contractible component of $\bar Y:=\tilde Y/P$
is $\bar B:=\tilde B/P$,
which is projected homeomorphically to $B$
by the covering projection $\bar p\co \bar X\to X$, where
$\bar X:=\tilde X/P$. 
Since $\bar X$  and $\bar B$ are aspherical CW-complexes,
and the inclusion
$\bar B\to\bar X$ is a $\pi_1$-isomorphism, this inclusion 
is a homotopy equivalence. Since $(\bar X,\bar B)$ is a CW-pair, 
it has the homotopy extension property,
so there is a deformation retraction $F_t\co\bar X\to\bar B$.
For the rest of the proof 
by a ``map'' we always mean a ``continuous map''.
To see that $\i\co C(Z,Y;f)\to C(Z,X;f)$
is a weak homotopy equivalence, it suffices 
to show that every map 
\[ H\co (D^{k+1},S^k)\to (C(Z,X;f),C(Z,Y;f))\] 
can be homotoped $rel\ S^k$ into 
$C(Z,Y;f)$~\cite[Lemma II.3.1]{Whi}.
 Since $C(Z,Y;f)$, $C(Z,X;f)$ are path-connected, we may assume
$k+1>0$. Fix a basepoint $\ast\in S^k$. 
Since $X$ is locally-compact, 
$H$ gives rise to a map $h\co D^{k+1}\times Z\to X$
that takes $S^k\times Z$ to $Y$ and whose restriction to
$\ast\times Z$ is a map in $C(Z,Y;f)=C(Z,B;f)$. 
Since $D^{k+1}$ is simply-connected,
$h$ lifts to $\bar h \co D^{k+1}\times Z\to\bar X$,
whose restriction to
$\ast\times Z$ is a map in $C(Z,\bar B;\bar f)$, where
$\bar f\co Z\to \bar B\subset\bar X$ be the unique
lift of $f$ to $\bar B$. 
If $k>0$, then $S^k\times Z$ is path-connected, and so
$\bar h(S^k\times Z)$ lies in $\bar B$.
Then $\bar p\circ F_t\circ\bar h$ is a homotopy 
$rel\ S^k\times Z$ from $h$ to a map with
image in $B$. 
If $k=0$, and $\ast =0\in \{0,1\}$, then 
$\bar h\co [0,1]\times Z\to\bar X$ is a homotopy between
the maps which we denote $\bar h_0$, $\bar h_1$ where 
$\bar h_0$ lies in $C(Z,\bar B;\bar f)$. 
The image of  $\bar h_1$ lies
in $\bar p^{-1}(B)$, because it projects to $B$. 
If the image of  $\bar h_1$ were in 
a contractible component of $p^{-1}(B)$, then
 $\bar h_1$ would be null-homotopic in $\bar X$, and so 
$\bar f$ would be null-homotopic in $\bar B$ 
which contradicts the assumption. 
Thus the image of  $\bar h_1$ lies in $\bar B$, and again
$\bar p\circ F_t\circ\bar h$ is a homotopy 
$rel\ S^k\times Z$ from $h$ to a map with
image in $B$. 
Thus $H$ is always homotopic $rel\ S^k$ to a map
with image in $C(Z,Y;f)$, as promised.
\end{proof}

\section{Hyperbolic manifolds with geodesic boundary}
\label{sec: geod bound}

The following result can be deduced from~\cite[Section 7]{Bow-rel}.
We find it worthwhile to spell out the details.

\begin{prop}\label{prop: geod bound}
If $M$ is a compact real hyperbolic
manifold with \textup(smooth\textup) totally geodesic boundary, then
$M$ satisfies Assumption~\ref{intro: assum}.
\end{prop}
\begin{proof}
The manifold $M$ is clearly locally convex, so
the universal cover $\tilde M$ of $M$ can be identified with
the convex complete codimension zero submanifold
of the real hyperbolic $n$-space $\bf{H^n}$
(see e.g. Corollary 1.3.7 and Proposition 1.4.2 in~\cite{CEG}).
The boundary components of $M$
are locally convex, hence they are
$\pi_1$-injectively embedded and aspherical, because 
they are totally geodesic and $\sec(M)\le 0$. 
Then the boundary components
of $\tilde M$ are simply-connected totally geodesic codimension 
one submanifolds of $\bf{H^n}$, thus each component of $\d\tilde M$
is a copy of $\bf{H^{n-1}}$. It follows that 
$\bf{H^n}\setminus\tilde M$ is the union of disjoint open $n$-dimensional
halfspaces in $\bf{H^n}$. We denote the halfspaces by $Q(p)$, $p\in\Pi$.
Clearly, $G$ acts on the set $Q$ of $Q(p)$'s, and the stabilizer of each
$Q(p)$ is isomorphic to the fundamental group of the corresponding
boundary component of $M$.
Bowditch proved in Proposition 7.12 and Lemma 7.13
of~\cite{Bow-rel} that given a quasidense locally finite
collection of uniformly convex subsets with bounded penetration
in a hyperbolic metric space, the nerve of the collection is a fine
hyperbolic graph. (The terms ``quasidense'', ``locally finite'', and
``uniformly convex'' are self-explanatory, while ``bounded penetration''
will be defined below).
Thus it is enough to show that $Q$  
satisfy the the above conditions. Each $Q(p)$ is convex, hence
uniformly convex. 
Since $M$ is compact, $Q$ is quasidense. Local finiteness of $Q$
easily follows from the fact that 
$G$-action on $\bf{H^n}$ is properly discontinuous. By definition
$Q$ has {\it bounded penetration} if for any $r\ge 0$ there exists
$D(r)$ such that for any distinct $Q(p), Q(q)\in Q$ the intersection of the
$r$-neighborhoods of $Q(p), Q(q)$ has diameter $\le D(r)$.
Arguing by contradiction assume that $Q$ does not have
bounded penetration. This means that for some $r$
there exist a sequence $x_i, y_i\in\bf{H^n}$ and $Q(p_i), Q(q_i)\in Q$
such that $x_i, y_i$ lie in the intersection of the
$r$-neighborhoods of $Q(p_i), Q(q_i)$, while the distance
$d(x_i,y_i)$ between $x_i$ and $y_i$ is $>i$. 
Given $z\in\bf{H^n}$, denote by $z(p)$ the orthogonal projection 
of $z$ to $Q(p)$. The geodesic segments
$[x_i(p_i), y_i(p_i)]$ and $[x_i(q_i), y_i(q_i)]$
have endpoints within $2r$, while their lengths are $\ge i-2r$. 
Because we are in a hyperbolic space
the segments become arbitrary close near the middle
for sufficiently large $i$, more precisely, their midpoints 
$m_{p_i}$, $m_{q_i}$ satisfy $d(m_{p_i},m_{q_i})\to 0$ 
as $i\to\infty$, and the same is true for any pair
of corresponding points on the segments that are within
say $i/4$ of $m_{p_i}$, $m_{q_i}$.
Since $M$ is compact, we can assume after passing to a subsequence that 
$m_{p_i}$ is nearly a constant sequence in $\bf{H^n}$. 
Then using that $Q$ is locally finite, 
we can pass to a subsequence to arrange that $Q(p_i), Q(q_i)$ 
are constant sequences. Thus $Q(p_i), Q(q_i)$ have to touch 
along the geodesic that is the pointed Hausdorff limit
of (either of) the segments. If two halfspaces in $\bf{H^n}$
touch along a geodesic, they intersect. So $Q(p_i), Q(q_i)$
intersect, which contradicts the fact that distinct
elements of $Q$ are disjoint.
\end{proof}

\begin{rmk} An analog of
Proposition~\ref{prop: geod bound} remains true (with a similar proof)
for complete finite volume Riemannian manifolds  
with pinched negative curvature and
compact totally geodesic boundary, namely, $\pi_1(M)$ is
hyperbolic relative to the fundamental groups of
boundary components and cusps of $M$. More generally, it should be
possible to prove the following relatively hyperbolic analog of 
Theorem 7.11 in~\cite{Bow-rel}: if $G$ is hyperbolic relative to 
$\mathcal P$,
and $\mathcal G$ is a conjugacy invariant family of infinite
quasiconvex subgroups such that any element of $\mathcal G$
equals to its normalizer, and the intersection of any two elements
of $\mathcal G$ is finite, then $G$ is hyperbolic relative to 
$\mathcal P\cup\mathcal G$. However, the focus of this paper 
is on manifolds, so I will not attempt to proof this
more general result here.
\end{rmk}

\section{Acknowledgements}

The author was partially supported by the NSF grants \# DMS-0503864,
DMS-0352576, and
is grateful to Ken Bromberg, Richard Canary, Michael Davis, Daniel Groves,
Ilya Kapovich, Vitali Kapovitch, Thilo Kuessner, John Millson,
Igor Mineyev, Denis Osin, and Mark Sapir for helpful 
discussions or communications.

\small
\bibliographystyle{amsalpha}
\bibliography{hp-r2}

\newcommand{\etalchar}[1]{$^{#1}$}
\def\cprime{$'$}
\providecommand{\bysame}{\leavevmode\hbox to3em{\hrulefill}\thinspace}
\providecommand{\MR}{\relax\ifhmode\unskip\space\fi MR }
\providecommand{\MRhref}[2]{%
  \href{http://www.ams.org/mathscinet-getitem?mr=#1}{#2}
}
\providecommand{\href}[2]{#2}
\begin{thebibliography}{ABC{\etalchar{+}}96}

\bibitem[AB05]{AB}
S.~B. Alexander and R.~L. Bishop, \emph{A cone splitting theorem for
  {A}lexandrov spaces}, Pacific J. Math. \textbf{218} (2005), no.~1, 1--15.

\bibitem[ABC{\etalchar{+}}96]{ABCKT}
J.~Amor{\'o}s, M.~Burger, K.~Corlette, D.~Kotschick, and D.~Toledo,
  \emph{Fundamental groups of compact {K}\"ahler manifolds}, Mathematical
  Surveys and Monographs, vol.~44, American Mathematical Society, Providence,
  RI, 1996.

\bibitem[BDM]{BerDruMos}
J.~Behrstock, C.~Dru{\c{t}}u, and L.~Mosher, \emph{{Thick metric spaces,
  relative hyperbolicity, and quasi-isometric rigidity}},
  arXiv:math.GT/0512592.

\bibitem[Bel02]{Bel-Top}
I.~Belegradek, \emph{On {M}ostow rigidity for variable negative curvature},
  Topology \textbf{41} (2002), no.~2, 341--361.

\bibitem[BF95]{BesFei}
M.~Bestvina and M.~Feighn, \emph{Stable actions of groups on real trees},
  Invent. Math. \textbf{121} (1995), no.~2, 287--321.

\bibitem[BH99]{BH}
M.~R. Bridson and A.~Haefliger, \emph{Metric spaces of non-positive curvature},
  Grundlehren der Mathematischen Wissenschaften [Fundamental Principles of
  Mathematical Sciences], vol. 319, Springer-Verlag, 1999.

\bibitem[BM99a]{Burger-Monod}
M.~Burger and N.~Monod, \emph{Bounded cohomology of lattices in higher rank
  {L}ie groups}, J. Eur. Math. Soc. (JEMS) \textbf{1} (1999), no.~2, 199--235.

\bibitem[BM99b]{Burger-Monod-err}
\bysame, \emph{Erratum: ``{B}ounded cohomology of lattices in higher rank {L}ie
  groups''}, J. Eur. Math. Soc. (JEMS) \textbf{1} (1999), no.~3, 338.

\bibitem[Bow]{Bow-rel}
B.~H. Bowditch, \emph{Relatively hyperbolic groups}, Southampton preprint,
  1999, www.maths.soton.ac.uk/staff/Bowditch/preprints.html.

\bibitem[BP92]{Ben-Pet}
R.~Benedetti and C.~Petronio, \emph{Lectures on hyperbolic geometry},
  Universitext, Springer-Verlag, 1992.

\bibitem[CD95]{CD}
R.~M. Charney and M.~W. Davis, \emph{Strict hyperbolization}, Topology
  \textbf{34} (1995), no.~2, 329--350.

\bibitem[CEG87]{CEG}
R.~D. Canary, D.~B.~A. Epstein, and P.~Green, \emph{Notes on notes of
  {T}hurston}, Analytical and geometric aspects of hyperbolic space
  (Coventry/Durham, 1984), London Math. Soc. Lecture Note Ser., vol. 111,
  Cambridge Univ. Press, Cambridge, 1987, pp.~3--92.

\bibitem[CM04]{CanMcC}
R.~D. Canary and D.~McCullough, \emph{Homotopy equivalences of 3-manifolds and
  deformation theory of {K}leinian groups}, Mem. Amer. Math. Soc. \textbf{172}
  (2004), no.~812, xii+218.

\bibitem[Dav02]{Davis-ex}
M.~W. Davis, \emph{Exotic aspherical manifolds}, Topology of high-dimensional
  manifolds, No. 1, 2 (Trieste, 2001), ICTP Lect. Notes, vol.~9, Abdus Salam
  Int. Cent. Theoret. Phys., Trieste, 2002, pp.~371--404.

\bibitem[Der05]{Deraux}
M.~Deraux, \emph{A negatively curved {K}\"ahler threefold not covered by the
  ball}, Invent. Math. \textbf{160} (2005), no.~3, 501--525.

\bibitem[DG05]{DelGro}
T.~Delzant and M.~Gromov, \emph{Cuts in {K}\"ahler groups}, Infinite groups:
  geometric, combinatorial and dynamical aspects, Progr. Math., vol. 248,
  Birkh\"auser, Basel, 2005, pp.~31--55.

\bibitem[DJ91]{DJ}
M.~W. Davis and T.~Januszkiewicz, \emph{Hyperbolization of polyhedra}, J.
  Differential Geom. \textbf{34} (1991), no.~2, 347--388.

\bibitem[DJW01]{DJW}
M.~W. Davis, T.~Januszkiewicz, and S.~Weinberger, \emph{Relative
  hyperbolization and aspherical bordisms: an addendum to ``{H}yperbolization
  of polyhedra'' [{J}.\ {D}ifferential {G}eom.\ {\bf 34} (1991), no.\ 2,
  347--388; {MR}1131435 (92h:57036)] by {D}avis and {J}anuszkiewicz}, J.
  Differential Geom. \textbf{58} (2001), no.~3, 535--541.

\bibitem[dlHV89]{dlHarVal}
P.~de~la Harpe and A.~Valette, \emph{La propri\'et\'e {$(T)$} de {K}azhdan pour
  les groupes localement compacts (avec un appendice de {M}arc {B}urger)},
  Ast\'erisque (1989), no.~175, 158, With an appendix by M. Burger.

\bibitem[DS]{DruSap-rips}
C.~Dru{\c{t}}u and M.~Sapir, \emph{{Groups acting on tree-graded spaces and
  splittings of relatively hyperbolic group}}, arXiv:math.GR/0601305.

\bibitem[DS05]{DOS}
C.~Dru{\c{t}}u and M.~Sapir, \emph{Tree-graded spaces and asymptotic cones of
  groups}, Topology \textbf{44} (2005), no.~5, 959--1058, With an appendix by
  D. Osin and Sapir.

\bibitem[Far98]{Far-rel}
B.~Farb, \emph{Relatively hyperbolic groups}, Geom. Funct. Anal. \textbf{8}
  (1998), no.~5, 810--840.

\bibitem[FJ89]{FarJon}
F.~T. Farrell and L.~E. Jones, \emph{Negatively curved manifolds with exotic
  smooth structures}, J. Amer. Math. Soc. \textbf{2} (1989), no.~4, 899--908.

\bibitem[FJO98]{FJO}
F.~T. Farrell, L.~E. Jones, and P.~Ontaneda, \emph{Hyperbolic manifolds with
  negatively curved exotic triangulations in dimensions greater than five}, J.
  Differential Geom. \textbf{48} (1998), no.~2, 319--322.

\bibitem[Fra04]{Fra}
S.~Francaviglia, \emph{Hyperbolic volume of representations of fundamental
  groups of cusped 3-manifolds}, Int. Math. Res. Not. (2004), no.~9, 425--459.

\bibitem[Fuj98]{Fuj}
K.~Fujiwara, \emph{The second bounded cohomology of a group acting on a
  {G}romov-hyperbolic space}, Proc. London Math. Soc. (3) \textbf{76} (1998),
  no.~1, 70--94.

\bibitem[Gol99]{Gol-nov}
B.~Goldfarb, \emph{Novikov conjectures and relative hyperbolicity}, Math.
  Scand. \textbf{85} (1999), no.~2, 169--183.

\bibitem[GPS88]{P-SGro}
M.~Gromov and I.~Piatetski-Shapiro, \emph{Nonarithmetic groups in {L}obachevsky
  spaces}, Inst. Hautes \'Etudes Sci. Publ. Math. (1988), no.~66, 93--103.

\bibitem[Gro82]{Grom-vol-bounded-coh}
M.~Gromov, \emph{Volume and bounded cohomology}, Inst. Hautes \'Etudes Sci.
  Publ. Math. (1982), no.~56, 5--99 (1983).

\bibitem[Gro87]{Gro-hgr}
\bysame, \emph{Hyperbolic groups}, Essays in group theory, Math. Sci. Res.
  Inst. Publ., vol.~8, Springer, 1987, pp.~75--263.

\bibitem[GS92]{GroSch}
M.~Gromov and R.~Schoen, \emph{Harmonic maps into singular spaces and
  {$p$}-adic superrigidity for lattices in groups of rank one}, Inst. Hautes
  \'Etudes Sci. Publ. Math. (1992), no.~76, 165--246.

\bibitem[GT87]{GroThu}
M.~Gromov and W.~Thurston, \emph{Pinching constants for hyperbolic manifolds},
  Invent. Math. \textbf{89} (1987), no.~1, 1--12.

\bibitem[Han81]{Han}
V.~L. Hansen, \emph{Spaces of maps into {E}ilenberg-{M}ac {L}ane spaces},
  Canad. J. Math. \textbf{33} (1981), no.~4, 782--785.

\bibitem[Hat02]{Hat}
A.~Hatcher, \emph{Algebraic topology}, Cambridge University Press, Cambridge,
  2002.

\bibitem[HR82]{HamRoy}
G.~C. Hamrick and D.~C. Royster, \emph{Flat {R}iemannian manifolds are
  boundaries}, Invent. Math. \textbf{66} (1982), no.~3, 405--413.

\bibitem[JY86]{Jos-Yau}
J.~Jost and S.-T. Yau, \emph{The strong rigidity of locally symmetric complex
  manifolds of rank one and finite volume}, Math. Ann. \textbf{275} (1986),
  no.~2, 291--304.

\bibitem[Kap01]{Kap-book}
M.~Kapovich, \emph{Hyperbolic manifolds and discrete groups}, Progress in
  Mathematics, vol. 183, Birkh\"auser Boston Inc., Boston, MA, 2001.

\bibitem[KS97]{KorSch}
N.~J. Korevaar and R.~M. Schoen, \emph{Global existence theorems for harmonic
  maps to non-locally compact spaces}, Comm. Anal. Geom. \textbf{5} (1997),
  no.~2, 333--387.

\bibitem[Kue]{Kue-mul}
T.~Kuessner, \emph{{Multicomplexes, bounded cohomology and additivity of
  simplicial volume}}, arXiv:math.AT/0109057.

\bibitem[LR00]{LonRei-eta}
D.~D. Long and A.~W. Reid, \emph{On the geometric boundaries of hyperbolic
  {$4$}-manifolds}, Geom. Topol. \textbf{4} (2000), 171--178.

\bibitem[LR02]{LonRei-orbi}
\bysame, \emph{All flat manifolds are cusps of hyperbolic orbifolds}, Algebr.
  Geom. Topol. \textbf{2} (2002), 285--296.

\bibitem[McR04]{McR}
D.~B. McReynolds, \emph{Peripheral separability and cusps of arithmetic
  hyperbolic orbifolds}, Algebr. Geom. Topol. \textbf{4} (2004), 721--755.

\bibitem[MS80]{MosSiu}
G.~D. Mostow and Y.-T. Siu, \emph{A compact {K}\"ahler surface of negative
  curvature not covered by the ball}, Ann. of Math. (2) \textbf{112} (1980),
  no.~2, 321--360.

\bibitem[MY]{MinYam}
I.~Mineyev and A.~Yaman, \emph{Relative hyperbolicity and bounded cohomology},
  preprint, 2006, www.math.uiuc.edu/~mineyev/.

\bibitem[Osi]{Osi-sc}
D.V. Osin, \emph{{Small cancellations over relatively hyperbolic groups and
  embedding theorems}}, arXiv:math.GR/0411039.

\bibitem[Osi06]{Osi-rel}
D.~V. Osin, \emph{Relatively hyperbolic groups: intrinsic geometry, algebraic
  properties, and algorithmic problems}, Mem. Amer. Math. Soc. \textbf{179}
  (2006), no.~843, vi+100.

\bibitem[Pau91]{Pau}
F.~Paulin, \emph{Outer automorphisms of hyperbolic groups and small actions on
  {${\bf R}$}-trees}, Arboreal group theory (Berkeley, CA, 1988), Math. Sci.
  Res. Inst. Publ., vol.~19, Springer, New York, 1991, pp.~331--343.

\bibitem[RS94]{RipSel}
E.~Rips and Z.~Sela, \emph{Structure and rigidity in hyperbolic groups. {I}},
  Geom. Funct. Anal. \textbf{4} (1994), no.~3, 337--371.

\bibitem[Spa81]{Spa}
E.~H. Spanier, \emph{Algebraic topology}, Springer-Verlag, New York, 1981,
  Corrected reprint.

\bibitem[Sun03]{Sun}
X.~Sun, \emph{Regularity of harmonic maps to trees}, Amer. J. Math.
  \textbf{125} (2003), no.~4, 737--771.

\bibitem[Swe01]{Swe}
E.~L. Swenson, \emph{Quasi-convex groups of isometries of negatively curved
  spaces}, Topology Appl. \textbf{110} (2001), no.~1, 119--129, Geometric
  topology and geometric group theory (Milwaukee, WI, 1997).

\bibitem[Szc02]{Szc-rel}
A.~Szczepa{\'n}ski, \emph{Examples of relatively hyperbolic groups}, Geom.
  Dedicata \textbf{93} (2002), 139--142.

\bibitem[Thu]{Thu-notes}
W.~P. Thurston, \emph{The geometry and topology of three-manifolds}, version
  2002, http://www.msri.org/publications/books/gt3m/.

\bibitem[Tuk94]{Tuk}
P.~Tukia, \emph{Convergence groups and {G}romov's metric hyperbolic spaces},
  New Zealand J. Math. \textbf{23} (1994), no.~2, 157--187.

\bibitem[Whi78]{Whi}
G.~W. Whitehead, \emph{Elements of homotopy theory}, Graduate Texts in
  Mathematics, vol.~61, Springer-Verlag, New York, 1978.

\bibitem[Yam97]{Yamag}
T.~Yamaguchi, \emph{Simplicial volumes of {A}lexandrov spaces}, Kyushu J. Math.
  \textbf{51} (1997), no.~2, 273--296.

\bibitem[Yam04]{Yaman}
A.~Yaman, \emph{A topological characterisation of relatively hyperbolic
  groups}, J. Reine Angew. Math. \textbf{566} (2004), 41--89.

\bibitem[Yeu91]{Yeu}
S.-K. Yeung, \emph{Compactification of {K}\"ahler manifolds with negative
  {R}icci curvature}, Invent. Math. \textbf{106} (1991), no.~1, 13--25.

\end{thebibliography}
\end{document}